\documentclass[letterpaper,onecolumn,12pt]{article}    
\usepackage{amssymb}
\usepackage{amsmath}
\usepackage{graphicx}  
\usepackage{epstopdf}
\usepackage{hyperref}
\usepackage{color}
\usepackage{bbold}
\usepackage{units}
\usepackage{a4wide}

\newtheorem{theorem}{Theorem}
\newtheorem{corollary}{Corollary}
\newtheorem{lemma}{Lemma}
\newtheorem{remark}{Remark}

\newtheorem{proposition}{Proposition}
\newtheorem{definition}{Definition}
\newtheorem{example}{Example}

\definecolor{changecolor}{rgb}{0.0, 0.0, 0.0}

\newcommand{\defeq}{\overset{\mathrm{def}}{=}}

\newcommand{\tr}{\intercal}

\newcommand{\proof}{\textbf{Proof. }}
\newcommand{\qed}{$\diamond$}

\title{\textbf{Intrinsic Separation Principles}}

\author{Boris Houska}
\date{\small ShanghaiTech University}

\begin{document}

\addtolength{\textheight}{-2cm}

\maketitle


\begin{abstract}
This paper is about output-feedback control problems for general linear systems in the presence of given state-, control-, disturbance-, and measurement error constraints. Because the traditional separation theorem in stochastic control is inapplicable to such constrained systems, a novel information-theoretic framework is proposed. It leads to an intrinsic separation principle that can be used to break the dual control problem for constrained linear systems into a meta-learning problem that minimizes an intrinsic information measure and a robust control problem that minimizes an extrinsic risk measure. The theoretical results in this paper can be applied in combination with modern polytopic computing methods in order to approximate a large class of dual control problems by finite-dimensional convex optimization problems.
\end{abstract}

\section{Introduction}
\label{sec:introduction}
The separation principle in stochastic control is a fundamental result in control theory~\cite{Joseph1961,Lindquist1973,Wonham1969}, closely related to the certainty-equivalence principle~\cite{Warter1981}. It states \mbox{that---under} certain \mbox{assumptions---the} problem of optimal control and state estimation can be decoupled.

\bigskip
\noindent
For general control systems, however, the separation theorem fails to hold. Thus, if one is interested in finding optimal output-feedback control laws for such systems, one needs to solve a rather complicated dual control problem~\cite{Feldbaum1961}. There are two cases where such dual- or output-feedback control problems are of interest:

\renewcommand{\theenumi}{\roman{enumi}}%
\begin{enumerate}
\addtolength{\itemsep}{2pt}

\item The first case is that we have an uncertain nonlinear system---in the easiest case, without state- and control constraints---for which the information content of future measurements depends on the control actions. In practice, this dependency can often be neglected, because, at least for small measurement errors and process noise, and under certain regularity assumptions, the separation theorem holds in a first order approximation~\cite{Stengel1994}. Nevertheless, there are some nonlinear systems that can only be stabilized if this dependency is taken into account~\cite{Filatov2004,Sehr2019}.

\item And, the second case is that we have an uncertain linear system with state- and control constraints. Here, the process noise and future measurement errors have to be taken into account if one wants to operate the system safely, for instance, by designing output-feedback laws that ensure constraint satisfaction for all possible uncertainty scenarios.
\end{enumerate}

\bigskip
\noindent
The current paper is about the second case. This focus is motivated by the recent trend towards the development of safe learning and control methods~\cite{Hewing2020,Zanon2021}.

\subsection{Literature Review}
Dual control problems have been introduced by Feldbaum in the early 1960s~\cite{Feldbaum1961}. Mature game-theoretic and stochastic methods for analyzing such dual- and output feedback control problems have, however, only been developed much later. They go back to the seminal work of N.N.~Krasovskii~\cite{Krasovskii1995,Krasovskii1964} and A.B.~Kurzhanskii,~\cite{Kurzhanski1972,Kurzhanski2004}. Note that these historical articles are complemented by modern set-theoretic control theory~\cite{Blanchini2003,Blanchini2015}. Specifically, in the context of constrained linear systems, set-theoretic notions of invariance under output feedback can be found in the work of D\'orea~\cite{Dorea2009,Dorea2021}, which focuses on the invariance of a single information set, and in the work of Artstein and Rakovi\'c~\cite{Artstein2011}, which focuses on the invariance of a collection of information sets. Moreover, a variety of set-theoretic output-feedback control methods for constrained linear systems have appeared in~\cite{Bemporad2000,Brunner2018,Efimov2022}. These have in common that they propose to take bounds on measurement errors into account upon designing a robust predictive controller. In this context, the work of Goulart and Kerrigan must be highlighted~\cite{Goulart2007}, who found a remarkably elegant way to optimize uncertainty-affine output feedback control laws for constrained linear systems. A general overview of methods for output-feedback and dual model predictive control (MPC) can be found in~\cite{Findeisen2003,Hovd2005,Mayne2009,Sehr2019}, and the reference therein.

\subsection{Contribution}

The three main contributions of this paper can be outlined as follows.

\bigskip
\noindent
\textbf{Meta Information Theory.} While traditional information theories are based on the assumption that one can learn from accessible data, models for predicting the evolution of an uncertain control system require a higher level of abstraction. Here, one needs a prediction structure that is capable of representing the set of all possible future information states of a dynamic learning process without having access to future measurement data. Note that a comprehensive and thorough discussion of this aspect can be found in the above mentioned article by Artstein and Rakovi\'c~\cite{Artstein2011}, in which notions of invariance under output-feedback for collections of information sets are introduced. Similar to their construction, the current article proposes a meta information theoretic framework that is based on a class of information set collections, too. A novel idea of the current article in this regard, however, is the introduction of intrinsic equivalence relations that can be used to categorize information sets with respect to their geometric properties. This leads to an algebraic-geometric definition of meta information spaces in which one can distinguish between extrinsic and intrinsic information measures. Here, intrinsic information about a system is needed to predict what we will know about its states, while extrinsic information is needed to predict and assess the risk that is associated to control decisions.

\bigskip
\noindent
\textbf{Intrinsic Separation Principle.} The central contribution of this paper is the introduction of the \textit{intrinsic separation principle}. It formalizes the fact that the intrinsic information content of a constrained linear system does not depend on the choice of the control law. An important consequence of this result is that a large class of dual receding horizon control problems can be solved by separating them into a meta learning problem that predicts intrinsic information and a robust control problem that minimizes extrinsic risk measures. Moreover, the intrinsic separation principle can be used to analyze the existence of solutions to dual control problems under certain assumptions on the continuity and monotonicity of the objective function of the dual control problem.

\bigskip
\noindent
\textbf{Polytopic Dual Control.} The theoretical results in this paper are used to develop practical methods to approximately solve dual control problems for linear systems with convex state- and control constraints as well as polytopic process noise and polytopic measurement error bounds. In order to appreciate the novelty of this approach, it needs to be recalled first that many existing robust output-feedback control methods, for instance the state-of-the-art output-feedback model predictive control methods in~\cite{Findeisen2003,Mayne2009}, are based on a set-theoretic or stochastic analysis of a coupled system-observer dynamics, where the control law depends on a state estimate. This is in contrast to the presented information theoretic approach to dual control, where control decisions are made based on the system's true information state rather than a state estimate. In fact, for the first time, this paper presents a polytopic dual control method that neither computes vector-valued state estimates nor introduces an affine observer structure. Instead, the discretization of the control law is based on optimizing a finite number of control inputs that are associated to so-called extreme polytopes. The shapes, sizes, and orientations of these extreme polytopes encode the system's intrinsic information while their convex hull encodes the system's extrinsic information. The result of this discretization is a finite dimensional convex optimization problem that approximates the original dual control problem.

\subsection{Overview}
The paper is structured as follows.
\begin{itemize}
\addtolength{\itemsep}{1pt}

\item Section~\ref{sec::InformationSpaces} reviews the main idea of set-theoretic learning and introduces related notation.

\item Section~\ref{sec::MetaLearning} establishes the technical foundation of this article. This includes the introduction of meta information spaces and a discussion of the difference between intrinsic and extrinsic information measures.

\item Section~\ref{sec::IntrinsicSP} introduces the intrinsic separation principle for constrained linear systems, see Theorem~\ref{thm::1}.

\item Section~\ref{sec::DualControl} discusses how to resolve dual control problems by intrinsic separation, see Theorem~\ref{thm::existence}.

\item Section~\ref{sec::Polytopes} presents methods for discretizing dual control problems using polytopic information set approximations. The main technical result is summarized in Theorem~\ref{thm::F}. A numerical case study is presented. And,

\item Section~\ref{sec::conclusions} summarizes the highlights of this paper.
\end{itemize}

\subsection{Notation}
Throughout this paper, $\mathbb K^n$ denotes the set of closed subsets of $\mathbb R^n$, while $\mathbb K_\mathrm{c}^n$ denotes the set of compact subsets of $\mathbb R^n$. It is equipped with the Hausdorff distance
\begin{align}
\notag
d_\mathrm{H}(X,Y) \defeq \max \left\{ \max_{x \in X} \min_{y \in Y} \| x-y \|, \max_{y \in Y} \min_{x \in X} \| x-y \| \right\}
\end{align}
for all $X,Y \in \mathbb K_\mathrm{c}^n$, where $\Vert \cdot \Vert: \mathbb R^n \to \mathbb R$ denotes a norm on~$\mathbb R^n$, such that $(\mathbb K_\mathrm{c}^n,d_\mathrm{H})$ is a metric space. This definition can be extended to $\mathbb K^n$ as follows: if the maxima in the above definition do not exist for \mbox{$X,Y \in \mathbb K^n$}, we set $d_\mathrm{H}(X,Y) = \infty$. The pair $(\mathbb K^n,d_\mathrm{H})$ is called an extended metric space.
Finally, the notation $\mathrm{cl}(\cdot)$ is used to denote the closure, assuming that it is clear from the context what the underlying metric distance function is. For instance, if $\mathfrak X \subseteq \mathbb K^{n}$ denotes a set of closed sets, $\mathrm{cl}(\mathfrak X)$ denotes the closure of $\mathfrak X$ in $(\mathbb K^n,d_\mathrm{H})$.

\section{Information Spaces}
\label{sec::InformationSpaces}
An information space $(\mathcal I,d,\sqcap)$ is a  space in which learning can take place. This means that $(\mathcal I,d)$ is an extended metric space that is equipped with a learning operator
$$\sqcap: \mathcal I \times \mathcal I \to \mathcal I \; ,$$
such that $(\mathcal I,\sqcap)$ is a semi-group. Among the most important examples for such spaces is the so-called set-theoretic information space, which is introduced below.

\subsection{Set-Theoretic Learning}
\label{sec::setLearning}
In the context of set-theoretic learning~\cite{Bertsekas1971,Blanchini2003,Witsenhausen1968}, $\mathcal I = \mathbb K^n$ denotes the set of closed subsets of the vector space~$\mathbb R^n$, while $d = d_\mathrm{H}$ denotes the (extended) Hausdorff distance. Here, the standard intersection operator takes the role of a learning operator, $$\sqcap = \cap \; ,$$
recalling that the intersection of closed sets is closed. The motivation behind this definition can be outlined as follows: let us assume that we currently know that a vector \mbox{$x \in \mathbb R^{n}$} is contained in a given set $X \in \mathbb K^n$. If we receive additional information, for instance, that the vector $x$ is also contained in the set $Y \in \mathbb K^n$, our posterior information is that $x$ is contained in the intersection of the sets $X$ and $Y$, which is denoted by $X \cap Y$.

\bigskip
\noindent
Note that the above set-theoretic framework is compatible with continuous functions. If $f: \mathbb R^n \to \mathbb R^m$ denotes such a continuous function, the notation
\[
\forall X \in \mathbb K^n, \qquad f(X) \ \defeq \ \{ f(x) \mid x \in X \}
\]
is used to denote its associated continuous image map. It maps closed sets in $\mathbb R^n$ to closed sets in $\mathbb R^m$. Similarly, for affine functions of the form $f(x) = A x + b$, the notation
\[
AX+b = \{ Ax+b \mid x \in X \}
\]
is used, where $A$ and $b$ are a matrix and a vector with compatible dimensions. And, finally, the Minkowski sum
\[
X+Y \ \defeq \ \{ x+y \mid x \in X, y \in Y \},
\]
is defined for all sets $X,Y \in \mathbb K^{n}$.

\begin{remark}
\label{rem::statLearning}
Set theoretic learning models can be augmented by probability measures in order to construct statistical information spaces~\cite{Doob1953}. In such a context, every element of $\mathcal I$ consists of a set $X$ and a probability distribution $\rho_X$ on~$X$. A corresponding metric is then constructed by using the Wasserstein distance~\cite{Villani2005}. Moreover, if $(X,\rho_X) \in \mathcal I$ and $(Y,\rho_Y) \in \mathcal I$ are two independent random variables, the learning operation
\[
(X,\rho_X) \sqcap (Y,\rho_Y) \ \defeq \ \left( \ X \cap Y , \ \rho_{XY} \ \right)
\]
has the form of a Bayesian learning update,
\[
\rho_{XY}(x) \ \defeq \ \frac{\rho_X(x)\rho_Y(x)}{\int_{X \cap Y} \rho_X(y)\rho_Y(y) \, \mathrm{d}y } \; .
\]
Thus, as much as the current paper focuses---for simplicity of presentation---on set-theoretic learning, most of the developments below can be generalized to statistical learning processes by augmenting the support sets with probability distributions or probability measures~\cite{Taylor1996,Wu2022}.
\end{remark}

\subsection{Expectation and Deviation}
\label{sec::deviation}
Expectation and deviation functions are among the most basic tools for analyzing learning processes~\cite{Doob1953}. The expectation function is defined by
\[
\forall X \in \mathbb K_\mathrm{c}^n, \qquad E(X) \ \defeq \ \int_X x \, \mathrm{d}x \; . 
\]
It is a continuous function on $\mathbb K_\mathrm{c}^n$ that satisfies
\[
E(AX+b) \ = \ A E(X) + b \; .
\]
For the special case that the compact set $X$ is augmented by its associated uniform probability distribution, as discussed in Remark~\ref{rem::statLearning}, the above definition of $E(X)$ corresponds to the traditional definition of expected value functions in statistics. Similarly, a deviation function $D: \mathbb K_\mathrm{c}^n \to \mathbb R$ is a continuous and radially unbounded function that satisfies
\begin{enumerate}
\addtolength{\itemsep}{2pt}
\item $D(X) \geq 0$,
\item $D(X) = 0$ if and only if $X = \{ E(X) \}$,
\item $D(X) = D(X-E(X))$, and
\item $D(X \cap Y) \leq D(X)$,
\end{enumerate}
for all $X,Y \in \mathbb K_\mathrm{c}^n$. While statistical learning models often use the variance of a random variable as a deviation measure, a more natural choice for $D$ in the context of set theoretic learning is given by the diameter,
\[
D(X) \ = \ \mathrm{diam}(X) \ \defeq \ \max_{x,y \in X} \ \| x - y \| \; .
\]
A long and creative list of other possible choices for $D$ can be found in~\cite{Rockafellar2013}.

\section{Meta Learning}
\label{sec::MetaLearning}
The above definition of information spaces assumes that information or data is accessible at the time at which learning operations take place. If one wishes to predict the future evolution of a learning process, however, one faces the problem that such data is not available yet. Therefore, this section proposes to introduce a meta information space in which one can represent the set of all possible posterior information states of a learning process without having access to its future data. Informally, one could say that a meta information space is an abstract space in which one can ``learn how to learn''.

\subsection{Information Ensembles}
The focus of this and the following sections is on the set-theoretic framework recalling that $\mathbb K_\mathrm{c}^n$ denotes the set of compact subsets of $\mathbb R^n$. A set $\mathfrak X \subseteq \mathbb K_\mathrm{c}^n$ is called an information ensemble of $\mathbb K_\mathrm{c}^n$ if
\begin{align}
\label{eq::iedef}
\forall Y \in \mathbb K_\mathrm{c}^n, \quad X \cap Y \in \mathfrak X
\end{align}
for all $X \in \mathfrak X$. Because $\varnothing = X \cap \varnothing \in \mathfrak X$, any information ensemble contains the empty set.

\bigskip
\noindent
\begin{proposition}
\label{prop::EnsembleClosure}
If $\mathfrak X \subseteq \mathbb K_\mathrm{c}^n$ is an information ensemble, then $\mathrm{cl}(\mathfrak X)$ is an information ensemble, too.
\end{proposition}

\bigskip
\noindent
\proof
Let $\mathfrak X$ be a given information ensemble and let $X_\infty \in \mathrm{cl}(\mathfrak X)$ be a given set in its closure. Then there exists a Cauchy sequence $X_1,X_2,\ldots \in \mathfrak X$ such that
\[
X_\infty \ \defeq \ \lim_{k \to \infty} X_k \ \in \ \mathrm{cl}(\mathfrak X) \; .
\]
Next, let $Y \in \mathbb K_\mathrm{c}^n$ be an arbitrary compact set. The case $X_\infty \cap Y = \varnothing$ is trivial, since $\varnothing \in \mathfrak X \subseteq \mathrm{cl}(\mathfrak X)$. Next, if $X_\infty \cap Y \neq \varnothing$, then there exists for every $\xi \in X_\infty \cap Y$ an associated sequence $z_1(\xi) \in X_1$, $z_2(\xi) \in X_2$, $\ldots \ $ with
\[
\lim_{k \to \infty} z_k(\xi) \ = \ \xi \; .
\]
This construction is such that the sets
\[
Z_k \ \defeq \ \mathrm{cl} \left( \ \{ z_k(\xi) \ \mid \ \xi \in X_\infty \cap Y \ \} \ \right)
\]
satisfy $Z_k = X_k \cap Z_k \in \mathfrak X$, since $\mathfrak X$ is an information ensemble. Consequently, it follows that
\[
 X_\infty \cap Y \ = \ \lim_{k \to \infty} Z_k \ \in \ \mathrm{cl}(\mathfrak X) \; .
\]
Thus, $\mathrm{cl}(\mathfrak X)$ is an information ensemble, as claimed by the statement of the proposition.
\qed

\bigskip
\noindent
Information ensembles can be used to construct information spaces, as pointed out below.

\bigskip
\noindent
\begin{proposition}
\label{prop::subsemigroup}
Let $\mathfrak X$ be an information ensemble of $\mathbb K_\mathrm{c}^n$. Then $(\mathfrak X,d_\mathrm{H},\cap)$ is an information space.
\end{proposition}

\bigskip
\noindent
\begin{proof}
Condition~\eqref{eq::iedef} implies that $X \cap Y \in \mathfrak X$ for all $X,Y \in \mathfrak X$. Thus, $(\mathfrak X,\cap)$ is a subsemigroup of $(\mathbb K_\mathrm{c}^n,\cap)$. Moreover, $d_\mathrm{H}$ defines a metric on $\mathfrak X$. Consequently, $(\mathfrak X,d_\mathrm{H},\cap)$ is an information space.\qed
\end{proof}

\bigskip
\noindent
\begin{remark}
The difference between information ensembles and more general set collections, as considered in~\cite{Artstein2011}, is that Property~\eqref{eq::iedef} is enforced. Note that this property makes a difference in the context of developing a coherent learning algebra: if~\eqref{eq::iedef} would not hold, $(\mathfrak X,\cap)$ would, in general, not be a subsemigroup of $(\mathbb K_\mathrm{c}^n,\cap)$.
\end{remark}

\subsection{Extreme Sets}
\label{sec::ExtemeSets}
A set $X \in \mathrm{cl}(\mathfrak X)$ of a given information ensemble $\mathfrak X \subseteq \mathbb K_\mathrm{c}^n$ is called an extreme set of $\mathfrak X$ if
\[
\forall Y \in \mathrm{cl}(\mathfrak X) \setminus \{ X \}, \qquad X \cap Y \neq X \; .
\]
The set of extreme sets of $\mathfrak X$ is denoted by $\partial \mathfrak X$. It is called the boundary of the information ensemble $\mathfrak X$. Clearly, we have $\partial X \subseteq \mathrm{cl}(\mathfrak X)$, but, in general, $\partial X$ is not an information ensemble. Instead, $\partial \mathfrak X$ can be interpreted as a minimal representation of the closure of $\mathfrak X$, because
\[
\mathrm{cl}(\mathfrak X) \ = \ \{ \ Y \in \mathbb K_\mathrm{c}^n \ \mid \ \exists X \in \partial \mathfrak X, \ Y \subseteq X \ \} \; .
\]
Reversely, the closure of $\mathfrak X$ can be interpreted as the smallest information ensemble that contains $\partial \mathfrak X$.

\subsection{Meta Information Spaces}
Let $\mathbb I^n$ denote the set of closed information ensembles of $\mathbb K_\mathrm{c}^n$; that is, the set of closed subsemigroups of $(\mathbb K_\mathrm{c}^n,\cap)$ that are closed under intersection with sets in $\mathbb K_\mathrm{c}^n$. Similarly, the notation $\mathbb I_\mathrm{c}^n$ will be used to denote the set of compact information ensembles of the information space $(\mathbb K_\mathrm{c}^n,d_\mathrm{H},\cap)$. Next, the meta learning operator \mbox{$\sqcap: \mathbb I^n \times \mathbb I^n \to \mathbb I^n$} is introduced by defining
\[
\mathfrak X \sqcap \mathfrak Y \ \defeq \ \{ \ X \cap Y \ \mid \ X \in \mathfrak X, \ Y \in \mathfrak Y \ \}
\]
for all $\mathfrak X,\mathfrak Y \in \mathbb I^n$. A corresponding metric distance function, $\Delta_H$ is given by
\begin{align}
& \Delta_\mathrm{H}(\mathfrak X, \mathfrak Y) \notag \\[0.16cm]
& \; \defeq \; \max \left\{ \max_{X \in \mathfrak X} \min_{Y \in \mathfrak Y} d_\mathrm{H}(X,Y), \max_{Y \in \mathfrak Y} \min_{X \in \mathfrak X} d_\mathrm{H}(X,Y) \right\} \notag
\end{align}
for all $\mathfrak X, \mathfrak Y \in \mathbb I_\mathrm{c}^n$ such that $(\mathbb I_\mathrm{c}^n,\Delta_\mathrm{H})$ is a metric space. Similar to the construction of the Hausdorff distance $d_\mathrm{H}$, the definition of $\Delta_\mathrm{H}$ can be extended to $\mathbb I^n$ by understanding the above definition in the extended value sense. The following proposition shows that the triple $(\mathbb I^n, \Delta_{\mathrm{H}}, \sqcap)$ is an information space. It is called the meta information space of $(\mathbb K_\mathrm{c}^n, d_{\mathrm{H}}, \cap)$. 

\bigskip
\noindent
\begin{proposition}
\label{prop::meta}
The triple $(\mathbb I^n, \Delta_{\mathrm{H}}, \sqcap)$ is an information space. It can itself be interpreted as a set-theoretic information space in the sense that we have
\begin{align}
\label{eq::cap}
\mathfrak X \sqcap \mathfrak Y \ = \ \mathfrak X \cap \mathfrak Y
\end{align}
for all $\mathfrak X,\mathfrak Y \in \mathbb I^n$.
\end{proposition}

\bigskip
\noindent
\begin{proof}
The proof of this proposition is divided into two parts: the first part shows that~\eqref{eq::cap} holds and the second part uses this result to conclude that $(\mathbb I^n, \Delta_{\mathrm{H}}, \sqcap)$ is an information space.

\bigskip
\noindent
\textit{Part I.} Let $\mathfrak X,\mathfrak Y \in \mathbb I^n$ be given information ensembles. For any $X \in \mathfrak X \cap \mathfrak Y$ the intersection relation
$$X \cap X = X \in \mathfrak X \cap \mathfrak Y$$ holds. But this implies that
\[
\mathfrak X \cap \mathfrak Y \ \subseteq \ \left\{ \ X \cap Y \ \middle| \ X \in \mathfrak X, \ Y \in \mathfrak Y \ \right\} \ = \ \mathfrak X \sqcap \mathfrak Y \; .
\]
In order to also establish the reverse inclusion, assume that $Z \in \mathfrak X \sqcap \mathfrak Y$ is a given set. It can be written in the form \mbox{$Z = X \cap Y$} with \mbox{$X \in \mathfrak X$} and \mbox{$Y \in \mathfrak Y$}. Clearly, we have \mbox{$Z \subseteq X$} and \mbox{$Z \subseteq Y$}. Moreover, we have \mbox{$Z \in \mathbb K_\mathrm{c}^n$}, since the intersection of compact sets is compact. Thus, since $\mathfrak X$ and $\mathfrak Y$ are information ensembles,~\eqref{eq::iedef} implies that \mbox{$Z \in \mathfrak X$} and \mbox{$Z \in \mathfrak Y$}. But this is the same as saying that \mbox{$Z \in \mathfrak X \cap \mathfrak Y$}, which implies $\mathfrak X \cap \mathfrak Y \supseteq \mathfrak X \sqcap \mathfrak Y$. Together with the above reverse inclusion, this yields~\eqref{eq::cap}.

\bigskip
\noindent
\textit{Part II.} Note that $(\mathbb I^n,\cap)$ is a semigroup, which follows from the definition of intersection operations. Moreover, $(\mathbb I^n,\Delta_\mathrm{H})$ is, by construction, an extended metric space. Thus, $(\mathbb I^n, \Delta_{\mathrm{H}}, \sqcap)$ is indeed an information space, as claimed by the statement of this proposition.\qed
\end{proof}

\bigskip
\noindent
\begin{corollary}
The triple $(\mathbb I_\mathrm{c}^n,\Delta_\mathrm{H},\sqcap)$ is also an information space. It can be interpreted as a sub-meta information space of $(\mathbb I^n,\Delta_\mathrm{H},\sqcap)$.
\end{corollary}

\bigskip
\noindent
\begin{proof}
The statement of this corollary follows immediately from the previous proposition, since the intersection of compact sets is compact; that is, $(\mathbb I_\mathrm{c}^n,\sqcap)$ is a subsemigroup of $(\mathbb I^n,\sqcap)$.\qed
\end{proof}

\bigskip
\noindent
The statement of the above proposition about the fact that $(\mathbb I^n,\Delta_\mathrm{H},\sqcap)$ can be interpreted as a set-theoretic information space can be further supported by observing that this space is naturally compatible with continuous functions, too. Throughout this paper, the notation
\[
f(\mathfrak X) \ \defeq \ \{ \ f(X) \ \mid \ X \in \mathfrak X \ \}
\]
is used for any $\mathfrak X \in \mathbb I^n$, recalling that $f(X)$ denotes the compact image set of a continuous function $f$ on a compact set $X \in \mathbb K_\mathrm{c}^n$. Due to this continuity assumption on $f$, closed information ensembles are mapped to closed information ensembles.

\subsection{Interpretation of Meta Learning Processes}
\label{sec::example}
Meta information spaces can be used to analyze the evolution of learning processes without having access to data. In order to discuss why this is so, a guiding example is introduced: let us consider a set-theoretic sensor, which returns at each time instance a compact information set $X \in \mathbb K_\mathrm{c}^1$ containing the scalar state $x$ of a physical system, $x \in X$. If the absolute value of the measurement error of the sensor is bounded by $1$, this means that $X \subseteq [a,a+2]$ for at least one lower bound $a \in \mathbb R$. The closed but unbounded information ensemble that is associated with such a sensor is given by
\begin{align}
\label{eq::sensorV}
\mathfrak V = \{ \ X \in \mathbb K_\mathrm{c}^1 \ \mid \ \exists a \in \mathbb R: \ X \subseteq [a,a+2] \ \}  \in \mathbb I^1 \; .
\end{align}
It can be interpreted as the set of all information sets that the sensor could return when taking a measurement.

\bigskip
\noindent
Next, in order to illustrate how an associated meta learning process can be modeled, one needs to assume that prior information about the physical state $x$ is available. For instance, if $x$ is known to satisfy $x \in [-3,3]$, this would mean that our prior is given by
\[
\mathfrak X = \{ \ X \in \mathbb K_\mathrm{c}^1 \ \mid \ X \subseteq [-3,3] \ \} \; .
\]
In such a situation, a meta learning process is---due to Proposition~\ref{prop::meta}---described by an update of the form
\[
\mathfrak X^+ \ = \ \mathfrak X \sqcap \mathfrak V \ = \ \mathfrak X \cap \mathfrak V \; ,
\]
where $\mathfrak X^+$ denotes the posterior,
\[
\mathfrak X^+ = \left\{ \ X \in \mathbb K_\mathrm{c}^1 \ \middle| \begin{array}{l}
\exists a \in \mathbb R: \\[0.16cm]
X \subseteq [ \max \{ a,-3\}, 2 + \min \{ 1, a \}  ]  
\end{array}
\right\} \; .
\]
It is computed without having access to any sensor data.

\subsection{Intrinsic Equivalence}
\label{sec::equivalence}
Equivalence relations can be used to categorize compact information sets with respect to their geometric properties. In the following, we focus on a particular equivalence relation. Namely, we consider two sets $X,Y \in \mathbb K_\mathrm{c}^n$ equivalent, writing $X \simeq Y$, if they have the same shape, size, and orientation. This means that
\[
X \simeq Y \qquad \Longleftrightarrow \qquad \exists a \in \mathbb R^{n}: \quad X + a = Y \; .
\]
The motivation for introducing this particular equivalence relation is that two information sets $X$ and $Y$ can be considered equally informative if they coincide after a translation.

\begin{definition}
\label{def::inequiv}
Two information ensembles $\mathfrak X,\mathfrak Y \subseteq \mathbb K_\mathrm{c}^{n}$ are called intrinsically equivalent, $\mathfrak X \sim \mathfrak Y$, if their quotient spaces coincide,
\begin{align}
\label{def::inequivDef}
(\mathfrak X/\simeq) \ = \ (\mathfrak Y/\simeq) \; .
\end{align}
\end{definition}
The intrinsic equivalence relation $\sim$ from the above definition is---as the name suggests---an equivalence relation. This follows from the fact that $\mathfrak X \sim \mathfrak Y$ if and only if
\begin{align}
\notag
\begin{array}{llllcl}
& \forall X \in \mathfrak X, & \exists a \in \mathbb R^n: & X + a & \in & \mathfrak Y \\[0.16cm]
\text{and} \quad & \forall Y \in \mathfrak Y, & \exists b \in \mathbb R^n: & Y + b & \in & \mathfrak X \; ,
\end{array}
\end{align}
which, in turn, follows after substituting the above definition of $\simeq$ in~\eqref{def::inequivDef}.

\bigskip
\noindent
\begin{proposition}
\label{prop::ClosureEquiv}
If $\mathfrak X,\mathfrak Y \subseteq \mathbb K_\mathrm{c}^n$ are intrinsically equivalent information ensembles, $\mathfrak X \sim \mathfrak Y$, their closures are intrinsically equivalent, too, 
\[
\mathrm{cl}(\mathfrak X) \ \sim \ \mathrm{cl}(\mathfrak Y) \; .
\]
\end{proposition}

\bigskip
\noindent
\proof
Proposition~\ref{prop::EnsembleClosure} ensures that the closures of $\mathfrak X$ and $\mathfrak Y$ are information ensembles, $\mathrm{cl}(\mathfrak X) \in \mathbb I^n$ and $\mathrm{cl}(\mathfrak Y) \in \mathbb I^n$. Next, there exists for every $X_\infty \in \mathrm{cl}(\mathfrak X)$ a convergent sequence of sets $X_1,X_2,\ldots \in \mathfrak X$ such that
\[
X_\infty = \lim_{k \to \infty} X_k \; .
\]
Moreover, since $\mathfrak X \sim \mathfrak Y$, there also exists a sequence $a_1,a_2,\ldots \in \mathbb R^n$ such that the sequence
\[
Y_k \ \defeq \ X_k + a_k \ \in \ \mathfrak Y 
\]
remains in $\mathfrak Y$. Because $\mathfrak X$ and $\mathfrak Y$ are compact, the sequence of offsets $a_k$ must be bounded. Thus, it has a convergent subsequence, $a_{j_1},a_{j_2},\ldots \in \mathbb R^{n}$, with limit
\[
a_\infty \ \defeq \ \lim_{k \to \infty} a_{j_k} \ \in \ \mathbb R^n \; .
\]
This construction is such that
\[
X_\infty + a_\infty \ = \ \lim_{k \to \infty} \ \{ X_{j_k} + a_{j_k} \} \ \in \ \mathrm{cl}(\mathfrak Y) \; .
\]
A completely analogous statement holds after replacing the roles of $\mathfrak X$ and $\mathfrak Y$. Consequently, the closures of $\mathfrak X$ and $\mathfrak Y$ are intrinsically equivalent, which corresponds to the statement of the proposition.
\qed

\subsection{Extrinsic versus Intrinsic Information}
\label{sec::intrinsic}
Throughout this paper, it will be important to distinguish between extrinsic and intrinsic information. Here, the extrinsic information of an information ensemble is encoded by the union of its elements, namely, the extrinsic information set. It describes present information. The extrinsic information content of an information ensemble can be quantified by extrinsic information measures:

\begin{definition}
\label{def::extrinsic}
An information measure $f: \mathbb I_\mathrm{c}^{n} \to \mathbb R$ is called extrinsic, if there exist a function $g: \mathbb K_\mathrm{c}^n \to \mathbb R$ with
\[
\forall \mathfrak X \in \mathbb I_\mathrm{c}^n, \qquad f(\mathfrak X) \ = \ g \left( \bigcup_{X \in \mathfrak X} X \right) \; .
\]
\end{definition}
In contrast to extrinsic information, the intrinsic information of an information ensemble $\mathfrak X$ is encoded by its quotient space, $\mathfrak X/\simeq$. It describes future information. In order to formalize this definition, it is helpful to introduce a shorthand for the meta quotient space
\[
\mathbb Q_\mathrm{c}^n \ \defeq \ \mathbb I_\mathrm{c}^n/\sim \; .
\]
In analogy to Definition~\ref{def::extrinsic}, the intrinsic information of an information ensemble can be quantified by intrinsic information measures:
\begin{definition}
An information measure $f: \mathbb I_\mathrm{c}^{n} \to \mathbb R$ is called intrinsic, if there exist a function $g: \mathbb Q_\mathrm{c}^n \to \mathbb R$ with
\[
\forall X \in \mathbb I_\mathrm{c}^n, \qquad f(\mathfrak X) \ = \ g(\mathfrak X/\simeq) \; .
\]
\end{definition}
In order to develop a stronger intuition about the difference between extrinsic and intrinsic information measures, it is helpful to extend the definitions of the expectation and deviation functions $E$ and $D$ from the original information space setting in Section~\ref{sec::deviation}. These original definitions can be lifted to the meta information space setting by introducing their associated extrinsic expectation $\mathfrak E$ and extrinsic deviation $\mathfrak D$, given by
\[
\mathfrak E(\mathfrak X) \ \defeq \ E\left( \bigcup_{X \in \mathfrak X} X \right) \quad \text{and} \quad \mathfrak D(\mathfrak X) \ \defeq \ D\left( \bigcup_{X \in \mathfrak X} X \right)
\]
for all $\mathfrak X \in \mathbb I_\mathrm{c}^n$. Note that $\mathfrak E$ and $\mathfrak D$ are continuous functions, which inherit the properties of $E$ and $D$. Namely, the relation
\[
\mathfrak E(A \mathfrak X + b ) = A \mathfrak E(\mathfrak X) + b
\]
holds. Similarly, $\mathfrak D$ satisfies all axioms of a deviation measure in the sense that
\begin{enumerate}
\addtolength{\itemsep}{2pt}
\item $\mathfrak D(\mathfrak X) \geq 0$,
\item $\mathfrak D(\mathfrak X) = 0$ if and only if $\mathfrak X = \{ \{ \mathfrak E( \mathfrak  X ) \} \}$,
\item $\mathfrak D(\mathfrak X) = \mathfrak D(\mathfrak X-\mathfrak E(\mathfrak X))$, and
\item $\mathfrak D(\mathfrak X \sqcap \mathfrak Y) \leq \mathfrak D(\mathfrak X)$,
\end{enumerate}
for all $\mathfrak X,\mathfrak Y \in \mathbb I_\mathrm{c}^n$. Note that such extrinsic deviation measures need to be distinguished carefully from intrinsic deviation measures. Here, a function $\mathfrak D^\circ: \mathbb I_\mathrm{c}^{n} \to \mathbb R$, is called an intrinsic deviation measure if it is a continuous and intrinsic function that satisfies
\begin{enumerate}
\addtolength{\itemsep}{2pt}
\item $\mathfrak D^\circ(\mathfrak X) \geq 0$,
\item $\mathfrak D^\circ(\mathfrak X) = 0$ if and only if $\mathfrak X \sim \{ \{ \mathfrak E( \mathfrak  X ) \} \}$,
\item $\mathfrak D^\circ(\mathfrak X) = \mathfrak D^\circ(\mathfrak X-\mathfrak E(\mathfrak X))$, and
\item $\mathfrak D^\circ(\mathfrak X \sqcap \mathfrak Y) \leq \mathfrak D^\circ(\mathfrak X)$,
\end{enumerate}
for all $\mathfrak X,\mathfrak Y \in \mathbb I_\mathrm{c}^n$. The second axiom is equivalent to requiring that $\mathfrak D^\circ$ is positive definite on the quotient space $\mathbb Q_\mathrm{c}^n$. In order to have a practical example in mind, we introduce the particular function
\begin{align}
\label{eq::Dcirc}
\forall \mathfrak X \in \mathbb I_\mathrm{c}^n, \qquad \mathfrak D_\infty^\circ(\mathfrak X) \ = \ \max_{X \in \mathfrak X} \ \max_{x,y \in X} \ \| x - y \| ,
\end{align}
which turns out to be an intrinsic information measure, as pointed out by the following lemma.

\bigskip
\noindent
\begin{lemma}
\label{lem::IntrinsicDeviation}
The function $\mathfrak D_\infty^\circ$, defined by~\eqref{eq::Dcirc}, is an intrinsic deviation measure on $\mathbb I_\mathrm{c}^n$.
\end{lemma}

\bigskip
\noindent
\proof
Let $\mathfrak X \in \mathbb I_\mathrm{c}^n$ be a given information ensemble and let $X^\star$ be a maximizer of~\eqref{eq::Dcirc}, such that
\begin{align}
\notag
\mathfrak D_\infty^\circ( \mathfrak X ) \ = \ \mathrm{diam}(X^\star) \ = \ \max_{x,y \in X^\star} \ \| x - y \| \; .
\end{align}
If $\mathfrak Y \in \mathbb I_\mathrm{c}^n$ is an intrinsically equivalent ensemble with $\mathfrak X \sim \mathfrak Y$, then there exists an offset vector $a^\star \in \mathbb R^{n}$ such that $X^\star+a^\star \in \mathfrak Y$. Thus, we have
\begin{eqnarray}
\mathfrak D_\infty^\circ(\mathfrak Y) & \ = \ & \max_{Y \in \mathfrak Y} \ \mathrm{diam}(Y) \ \geq \ \mathrm{diam}(X^\star + a^\star) \notag \\[0.16cm]
& =& \mathrm{diam}(X^\star + a - E(X^\star + a) ) \notag \\[0.16cm]
&=& \mathrm{diam}(X^\star - E(X^\star)) \notag \\[0.16cm]
&=& \mathrm{diam}(X^\star) = \mathfrak D_\infty^\circ(\mathfrak X) \; , \notag
\end{eqnarray}
where the equations in the second, third, and fourth line follow by using the axioms of $D$ and $E$ from Section~\ref{sec::deviation}. The corresponding reverse inequality follows by using an analogous argument exchanging the roles of $\mathfrak X$ and $\mathfrak Y$. Thus, we have $\mathfrak D_\infty^\circ(\mathfrak X) = \mathfrak D_\infty^\circ(\mathfrak Y)$. This shows that $\mathfrak D_\infty^\circ$ is an intrinsic information measure. The remaining required properties of $\mathfrak D_\infty^\circ$ are directly inherited from
 the diameter function, recalling that the diameter is a continuous deviation function that satisfies the corresponding axioms from Section~\ref{sec::deviation}. This yields the statement of the lemma. \qed

\begin{example}
Let us revisit the tutorial example from Section~\ref{sec::example}, where we had considered the case that
\begin{align}
\mathfrak X &= \{ \ X \in \mathbb K_\mathrm{c}^1 \ \mid \ X \subseteq [-3,3] \ \} \qquad \text{and} \notag \\[0.16cm]
\mathfrak X^+ &= \left\{ \ X \in \mathbb K_\mathrm{c}^1 \ \middle| \begin{array}{l}
\exists a \in \mathbb R: \\
X \subseteq [ \max \{ a,-3\}, 2 + \min \{ 1, a \}  ]  
\end{array}
\right\} \notag
\end{align}
denote the prior and posterior of a data-free meta learning process. If we set $D(X) = \mathrm{diam}(X)$ and define $\mathfrak D$ and $\mathfrak D_\infty^\circ$ as above, then
\[
\mathfrak D( \mathfrak X ) \ = \ \mathfrak D( \mathfrak X^+ ) \ = \ 6 \; . 
\]
An interpretation of this equation can be formulated as follows: since our meta learning process is not based on actual data, the extrinsic information content of the prior $\mathfrak X$ and the posterior $\mathfrak X^+$ must be the same, which implies that their extrinsic deviations must coincide. This is in contrast to the intrinsic deviation measure,
\[
\mathfrak D_\infty^\circ( \mathfrak X ) \ = \ 6 \ > \ 2 \ = \ \mathfrak D_\infty^\circ( \mathfrak X^+ ),
\]
which predicts that no matter what our next measurement will be, the diameter of our posterior information set will be at most $2$.
\end{example}

\section{Intrinsic Separation Principle}
\label{sec::IntrinsicSP}
The goal of this section is to formulate an intrinsic separation principle for constrained linear systems.

\subsection{Constrained Linear Systems}

The following considerations concern uncertain linear discrete-time control systems of the form
\begin{eqnarray}
\label{eq::sys}
\begin{array}{rcl}
x_{k+1} &=& A x_k + B u_k + w_k \\[0.16cm]
\eta_k &=& C x_k + v_k \; .
\end{array}
\end{eqnarray}
Here, $x_k \in \mathbb R^{n}$ denotes the state, $u_k\in \mathbb U$ the control, $w_k \in \mathbb W$ the disturbance, \mbox{$\eta_k \in \mathbb R^{n_v}$} the measurement, and $v_k \in \mathbb V$ the measurement error at time $k \in \mathbb Z$. The system matrices $A$, $B$, and $C$ as well as the state, control, disturbance, and measurement error constraints sets, $\mathbb X \in \mathbb K^{n}$, $\mathbb U \in \mathbb K_\mathrm{c}^{n_u}$, $\mathbb W \in \mathbb K_\mathrm{c}^{n}$, and $\mathbb V \in \mathbb K_\mathrm{c}^{n_v}$, are assumed to be given.

\subsection{Information Tubes}
The sensor that measures the outputs $C x_k$ of~\eqref{eq::sys} can be represented by the information ensemble
\begin{align}
\label{eq::V}
\mathfrak V \ \defeq \ \left\{ \, X \in \mathbb K_\mathrm{c}^{n} \ \middle|  \ \exists \eta \in \mathbb R^{n_v} \hspace{-0.05cm} : \eta -C X \subseteq \mathbb V \, \right\} .
\end{align}
Since $\mathbb V$ is compact, $\mathfrak V$ is closed but potentially unbounded, $\mathfrak V \in \mathbb I^{n}$. If $\mathfrak X \in \mathbb I^{n}$ denotes a prior information ensemble of the state of~\eqref{eq::sys} an associated posterior is given by $\mathfrak X \sqcap \mathfrak V$. This motivates to introduce the function
\[
F( \mathfrak X, \mu ) \ \defeq \ \left\{ X^+ \in \mathbb K_\mathrm{c}^{n} \middle|
\begin{array}{l}
\exists X \in \mathfrak X \sqcap \mathfrak V: \\[0.1cm]
X^+ \subseteq A X + B \mu(X) + \mathbb W
\end{array}
\right\},
\]
which is defined for all $\mathfrak X \in \mathbb I^n$ and all control laws \mbox{$\mu: \mathbb I^n \to \mathbb U$} that map the system's posterior information state to a feasible control input. Let $\mathcal U$ denote the set of all such maps from $\mathbb I^{n}$ to $\mathbb U$. It is equipped with the supremum norm,
\[
\Vert \mu \Vert \ \defeq \ \sup_{X \in \mathbb K^n} \ \left\| \mu(X) \right\| \; ,
\]
such that $(\mathcal U,\Vert \cdot \Vert)$ is a Banach space. As $\mu \in \mathcal U$ is potentially discontinuous, $F(\mathfrak X,\mu)$ is not necessarily closed. Instead, the following statement holds.

\bigskip
\noindent
\begin{proposition}
If $\mathfrak X$, $\mathbb U$, $\mathbb V$, and $\mathbb W$ are closed, then the closure of the set $F(\mathfrak X,\mu)$ is for every given $\mu \in \mathcal U$ a closed information ensemble,
\[
\overline F(\mathfrak X,\mu) \ \defeq \ \mathrm{cl}( \ F(\mathfrak X,\mu) \ ) \in \mathbb I^n \; .
\]
\end{proposition}

\bigskip
\noindent
\proof
The statement of this proposition follows from Proposition~\ref{prop::EnsembleClosure} and the above definition of~$F$.
\qed

\bigskip
\noindent
The functions $F$ and $\overline F$ are the basis for the following definitions.

\bigskip
\noindent
\begin{definition}
An information ensemble $\mathfrak X_\mathrm{s} \in \mathbb I^n$ is called control invariant~\eqref{eq::sys} if there exists a \mbox{$\mu_\mathrm{s} \in \mathcal U$} such that
\[
\mathfrak X_\mathrm{s} \supseteq F(\mathfrak X_\mathrm{s}, \mu_\mathrm{s}) \; .
\]
\end{definition}

\bigskip
\noindent
\begin{definition}
A sequence $\mathfrak X_0,\mathfrak X_1,\ldots \in \mathbb I^n$ of information ensembles is called an information tube for~\eqref{eq::sys} if there exists a sequence $\mu_0,\mu_1,\ldots \in \mathcal U$ such that
\[
\forall k \in \mathbb N, \quad \mathfrak X_{k+1} \supseteq F(\mathfrak X_k, \mu_k) \; .
\]
\end{definition}

\bigskip
\noindent
\begin{definition}
An information tube $\mathcal X_0,\mathcal X_1,\ldots \in \mathbb I^n$ is called tight if it satisfies
\[
\forall k \in \mathbb N, \qquad \mathfrak X_{k+1} = \overline F( \mathfrak X_k, \mu_k )
\]
for at least one control policy sequence $\mu_k \in \mathcal U$.
\end{definition}

\subsection{Intrinsic Separation}
The following theorem establishes the fact that the intrinsic equivalence class of tight information tubes does not depend on the control policy sequence.

\bigskip
\noindent
\begin{theorem}
\label{thm::1}
Let $\mathfrak X_0,\mathfrak X_1, \ldots \in \mathbb I_\mathrm{c}^n$ and $\mathfrak Y_0, \mathfrak Y_1, \ldots \in \mathbb I_\mathrm{c}^n$ be tight information tubes with compact elements. If the initial information ensembles are intrinsically equivalent, $\mathfrak X_0 \sim \mathfrak Y_0$, then all information ensembles are intrinsically equivalent; that is, $\mathfrak X_k \sim \mathfrak Y_k$ for all $k \in \mathbb N$.
\end{theorem}

\bigskip
\noindent
\proof
Because $\mathfrak X$ and $\mathfrak Y$ are tight information tubes, there exist control policies \mbox{$\mu_k: \mathfrak X_k \cap \mathfrak V \to \mathbb U$} and \mbox{$\nu_k: \mathfrak Y_k \cap \mathfrak V \to \mathbb U$} such that
\begin{eqnarray}
\label{eq::XY}
\mathfrak X_{k+1} = \overline F(\mathfrak X_k, \mu_k ) \quad \text{and} \quad \mathfrak Y_{k+1} = \overline F(\mathfrak Y_k, \nu_k )
\end{eqnarray}
for all $k \in \mathbb N$. Next, the statement of the theorem can be proven by induction over $k$: since we assume $\mathfrak X_0 \sim \mathfrak Y_0$, this assumption can be used directly as induction start. Next, if $\mathfrak X_k \sim \mathfrak Y_k$, there exists for every $X_k \in \mathfrak X_k \cap \mathfrak V$ an offset vector $a_k \in \mathbb R^n$ such that $Y_k = X_k + a_k \in \mathfrak Y_k$. Because $\mathfrak V$ satisfies
\[
\forall a \in \mathbb R^{n}, \ \forall V \in \mathfrak V, \qquad V + a \in \mathfrak V,
\]
it follows that $Y_k = X_k + a_k \in \mathfrak Y_k \cap \mathfrak V$. Consequently, a relation of the form
\begin{eqnarray}
A X_k + B \mu_k(X_k) + \mathbb W &=& A Y_k + (B\mu_k(X_k) - Aa_k) + \mathbb W \notag \\[0.16cm]
&=& A Y_k + B \nu_k(Y_k) + \mathbb W - a_{k+1}, \notag
\end{eqnarray}
can be established, where the next offset vector, $a_{k+1}$, is given by
\[
a_{k+1} \ \defeq \ A x_k + B \nu_k(Y_k) - B \mu_k(X_k) \ \in \ \mathbb R^n \; .
\]
Note that a completely symmetric relation holds after exchanging the roles of $\mathfrak X_k$ and $\mathfrak Y_k$. In summary, it follows that an implication of the form
\[
\mathfrak X_k \sim \mathfrak Y_k \qquad \Longrightarrow \qquad F(\mathfrak X_k,\mu_k) \sim F(\mathfrak Y_k,\nu_k)
\]
holds. An application of Proposition~\ref{prop::ClosureEquiv} to the latter equivalence relation yields the desired induction step. This completes the proof of the theorem.\qed

\bigskip
\noindent
The above theorem allows us to formulate an intrinsic separation principle. Namely, Theorem~\ref{thm::1} implies that the predicted future information content of a tight information tube does not depend on the choice of the control policy sequence with which it is generated. In particular, the tight information tubes from~\eqref{eq::XY} satisfy
\[
\forall k \in \mathbb N, \qquad \mathfrak D^\circ(\mathfrak X_k) = \mathfrak D^\circ(\mathfrak Y_k)
\]
for any intrinsic information measure $\mathfrak D^\circ$. Note that this property is independent of the choice of the control policy sequences $\mu_k$ and $\nu_k$ that are used to generate these tubes.

\subsection{Control Invariance}
As mentioned in the introduction, different notions of invariance under output-feedback control have been analyzed by various authors~\cite{Dorea2021,Artstein2011}. This section briefly discusses how a similar result can be recovered by using the proposed meta learning based framework. For this aim, we assume that
\renewcommand{\theenumi}{\arabic{enumi}}%
\begin{enumerate}
\addtolength{\itemsep}{2pt}

\item the sets $\mathbb V \in \mathbb K_\mathrm{c}^{n_v}$ and $\mathbb W \in \mathbb K_\mathrm{c}^{n_w}$ are compact,

\item the set \mbox{$\mathbb U \in \mathbb K^{n_u}$} is closed and convex,

\item the pair $(A,C)$ is observable, and

\item $(A,B,\mathbb U, \mathbb W)$ admits a robust control invariant set. 
\end{enumerate}

\bigskip
\noindent
The first two assumptions are standard. The third assumption on the observability of $(A,C)$ could also be replaced by a weaker detectability condition. However, since one can always use a Kalman decomposition to analyze the system's invariant subspaces separately~\cite{Kalman1962}, it is sufficient to focus on observable systems. And, finally, the fourth assumption is equivalent to requiring the existence of a state-feedback law $\overline \mu: \mathbb R^n \to \mathbb U$ and a set $\overline X \in \mathbb K_\mathrm{c}^{n}$ such that
$$\forall x \in \overline X, \ \forall w \in \mathbb W, \qquad Ax+B \overline \mu(x) +w \in \overline X \, ,$$
which is clearly necessary: if we cannot even keep the system inside a bounded region by relying on exact state measurements, there is no hope that we can do so without such exact data.  

\bigskip
\noindent
\begin{lemma}
\label{lem::InvariantSet}
If the above four assumptions hold,~\eqref{eq::sys} admits a compact control invariant information ensemble.
\end{lemma}

\bigskip
\noindent
\proof The proof of this lemma is divided into two parts, which aim at constructing an information tube that converges to a control invariant information ensemble.

\bigskip
\noindent
\textit{Part I.} The goal of the first part is to show, by induction over $k$, that the recursion
\[
\forall k \in \mathbb N, \qquad \mathfrak X_{k+1}^\circ \ \defeq \ A ( \mathfrak X_k^\circ \cap \mathfrak V ), \quad \mathfrak X_0^\circ \ \defeq \ \mathbb K_\mathrm{c}^n
\]
is set monotonous. Since $\mathfrak X_0^\circ = \mathbb K_\mathrm{c}^n$, $\mathfrak X_1^\circ \subseteq \mathfrak X_0^\circ$ holds. This is the induction start. Next, if $\mathfrak X_{k+1}^\circ \subseteq \mathfrak X_k^\circ$ holds for a given integer $k \geq 0$, it follows that
\begin{eqnarray}
\label{eq::10}
\mathfrak X_{k+2}^\circ \ = \ A(\mathfrak X_{k+1}^\circ \cap \mathfrak V ) \ \subseteq \ A(\mathfrak X_k^\circ \cap \mathfrak V ) = \mathfrak X_{k+1}^\circ \; ,
\end{eqnarray}
where the inclusion in the middle follows directly by substituting the induction assumption. In summary, the monotonicity relation $\mathfrak X_{k+1}^\circ \subseteq \mathfrak X_k^\circ$ holds for all $k \in \mathbb N$.

\bigskip
\noindent
\textit{Part II.} The goal of the second part is to show that the sequence
\begin{eqnarray}
\label{eq::auxXk}
\mathfrak X_k \, \defeq \, \left\{ \ X - E(X) + \overline x \ \middle| \ X \in \mathfrak X_k^\circ, \ \overline x \in \mathrm{cvx}(\overline X) \ \right\} ,
\end{eqnarray}
converges to an invariant information ensemble. Here, $\mathrm{cvx}(\overline X)$ denotes the convex hull of the robust control invariant set $\overline X$. Because we assume that $\mathbb U$ is convex, $\mathrm{cvx}(\overline X)$ is robust control invariant, too. This means that there exists a $\overline \mu: \mathbb R^n \to \mathbb U$ such that
\[
\forall x \in \mathrm{cvx}(\overline X), \forall w \in \mathbb W, \quad Ax + B \overline \mu(x)+w \in \mathrm{cvx}(\overline X) \; .
\]
Since $E$ satisfies $E( X ) \in \mathrm{cvx}(X)$ for all $X \in \mathbb K_\mathrm{c}^n$,~\eqref{eq::auxXk} and the definitions of $\mathfrak X_k$ and $\mathfrak V$ imply that
\begin{align}
\begin{array}{llrl}
& \forall X \in \mathfrak X_k \cap \mathfrak V, & \qquad E(X) & \in \mathrm{cvx}(\overline X), \\
& \forall X \in \mathfrak X_k, & \qquad X - E(X) & \in \mathfrak X_k^\circ \\
\text{and} \qquad & \forall X \in \mathfrak V, & \qquad X - E(X) & \in \mathfrak V
\end{array} \notag
\end{align}
for all $k \in \mathbb N$. Thus, the state estimation based auxiliary feedback law
\begin{align}
\label{eq::StateFeedback}
\forall X \in \mathbb K_\mathrm{c}^n, \qquad \mu(X) \ \defeq \ \overline \mu(E(X))
\end{align}
ensures that the recursive feasibility condition
\begin{eqnarray}
\begin{array}{l}
A X + B \mu(X) + \mathbb W \\
= \ A(X-E(X)) + \underbrace{A E(X) + B \overline \mu(E(X)) + \mathbb W}_{\subseteq \, \mathrm{cvx}(\overline X)} \in \mathfrak X_{k+1}
\end{array} \notag
\end{eqnarray}
holds for all $X \in \mathfrak X_k \cap \mathfrak V$. Consequently, the auxiliary sequence $\mathfrak X_k$ is a monotonous information tube,
\[
\forall k \in \mathbb N, \qquad \mathfrak X_k \ \supseteq \ \mathfrak X_{k+1} \ \supseteq \ F(\mathfrak X_k,\mu) \; ,
\]
where monotonicity follows from~\eqref{eq::10} and the considerations from Part I.
Moreover, since $(A,C)$ is observable, $\mathfrak X_k$ is compact for all $k \geq n-1$. In summary, $\mathfrak X_k$ is a monotonously decreasing sequence of information ensembles, which---due to the monotone convergence theorem---converges to a compact control invariant information ensemble,
\[
\mathfrak X_\infty \ = \ \lim_{k \to \infty} \ \mathfrak X_k \ \in \ \mathbb I_\mathrm{c}^{n_x} \quad \ \text{and} \quad \ F(\mathfrak X_\infty,\mu) \ \subseteq \ \mathfrak X_\infty \; .
\]
This corresponds to the statement of the lemma.
\qed

\bigskip
\noindent
\begin{remark}
\label{rem::InvariantSet}
The purpose of Lemma~\ref{lem::InvariantSet} is to elaborate on the relation between control invariant information ensembles and existing notions in linear control \mbox{theory---such} as observability and robust stabilizability. Lemma~\ref{lem::InvariantSet} does, however, not make statements about feasibility: the state constraint set $\mathbb X$ is not taken into account. Moreover, the construction of the feedback law $\mu$ in~\eqref{eq::StateFeedback} is based on the vector-valued state estimate $E(X)$ rather than the information state $X$, which is, in general, sub-optimal. Note that these problems regarding feasibility and optimality are resolved in the following section by introducing optimal dual control laws.
\end{remark}

\section{Dual Control}
\label{sec::DualControl}
This section is about dual control problems for constrained linear systems. It is discussed under which assumptions such problems can be separated into a meta learning and a robust control problem.

\subsection{Receding Horizon Control}
\label{sec::RHC}
Dual control problems can be implemented in analogy to traditional model predictive control (MPC) methods. Here, one solves the online optimization problem
\begin{eqnarray}
J(X_0) \ = \ &\; \; \inf_{\mathfrak X,\mu} \; \; & \sum_{k=0}^{N-1} L(\mathfrak X_k,\mu_k) + M(\mathfrak X_N)
\notag \\[0.2cm]
\label{eq::dcp}
&\mathrm{s.t.}& \left\{
\begin{array}{l}
\forall k \in \{ 0, 1, \ldots, N-1 \}, \\[0.1cm]
F(\mathfrak X_k,\mu_k) \subseteq \mathfrak X_{k+1}, \ \ X_0 \in \mathfrak X_0 \\[0.1cm]
\mu_k \in \mathcal U, \\
\forall X_k \in \mathfrak X_k,  \ \ X_k \subseteq \mathbb X 
\end{array}
\right.
\end{eqnarray}
on a finite time horizon $\{ 0,1,\ldots, N \}$, where $0$ is the current time. The optimization variables are the feedback policies $\mu_0,\mu_1,\ldots,\mu_{N-1} \in \mathcal U$ and their associated information tube, $\mathfrak X_0,\mathfrak X_1,\ldots,\mathfrak X_N \in \mathbb I_\mathrm{c}^{n}$. In the most general setting, the stage and terminal cost functions,
\[
L: \mathbb I_\mathrm{c}^{n} \times \mathcal U \to \mathbb R \qquad \text{and} \qquad M: \mathbb I_\mathrm{c}^{n} \to \mathbb R,
\]
are assumed to be lower semi-continuous, although some of the analysis results below will be based on stronger assumptions. We recall that $\mathbb X$ denotes the closed state constraint set. The parameter $X_0 \in \mathbb I_\mathrm{c}^{n}$ corresponds to the current information set. It is updated twice per sampling time by repeating the following steps online:

\begin{enumerate}
\addtolength{\itemsep}{2pt}

\item[i)] Wait for the next measurement $\eta$.
\item[ii)] Update the information set,
$$X_0 \leftarrow X_0 \cap \{ x \in \mathbb R^n \mid \eta - Cx \in \mathbb V \} \; .$$
\item[iii)] Solve~\eqref{eq::dcp} and denote the first element of the optimal feedback sequence by $\mu_0^\star \in \mathcal U$.

\item[iv)] Send $u^\star = \mu_0^\star(X_0)$ to the real process.

\item[v)] Propagate the information set,
$$X_0 \leftarrow A X_0 + B u^\star + \mathbb W \; .$$

\item[vi)] Set the current time to $0$ and continue with Step~i). 

\end{enumerate}

\bigskip
\noindent
Note that Step~iii) assumes that the ``$\inf$'' operator in~\eqref{eq::dcp} can be replaced by a ``$\min$'' and that an associated optimal feedback policy exists. Conditions under which this can be guaranteed are discussed in Section~\ref{sec::existence}.

\subsection{Objectives}

Tube model predictive control formulations~\cite{Rawlings2017,Villanueva2020,Villanueva2022} use risk measures as stage cost functions. In principle, any lower semi-continuous function of the form
\[
R: \mathbb K_\mathrm{c}^n \to \mathbb R \cup \{ \infty \},
\]
can be regarded as such a risk measure, although one would usually require that the monotonicity condition
\begin{align}
\label{eq::risk}
X \ \subseteq \ Y \qquad \Longrightarrow \qquad R(X) \ \leq \ R(Y) 
\end{align}
holds for all $X,Y \in \mathbb K_\mathrm{c}^n$. Similarly, this paper proposes to call $\mathfrak R: \mathbb I_\mathrm{c}^n \to \mathbb R \cup \{ \infty \}$ an extrinsic risk measure if
\begin{align}
\label{eq::metarisk}
\mathfrak R(\mathfrak X) \ = \ R \left( \ \bigcup_{X \in \mathfrak X} X \ \right)
\end{align}
for a lower semi-continuous function $R$ that satisfies~\eqref{eq::risk}.
\begin{remark}
Problem~\eqref{eq::dcp} enforces state constraints explicitly. Alternatively, one can move them to the objective by introducing the indicator function $I_{\mathbb X}$ of the state constraint set $\mathbb X$. Because we have
\[
\left( \ \forall X_k \in \mathfrak X_k,  \ \ X_k \subseteq \mathbb X \ \right) \quad \Longleftrightarrow \quad I_{\mathbb X} \left( \bigcup_{X \in \mathfrak X_k} X \right) < \infty,
\]
enforcing state constraints is equivalent to adding an extrinsic risk measure to the stage cost; here with $R = I_{\mathbb X}$. 
\end{remark}

\bigskip
\noindent
By using the language of this paper, the traditional objective of dual control~\cite{Feldbaum1961} is to tradeoff between extrinsic risk and intrinsic deviation. This motivates to consider stage cost functions of the form
\begin{align}
\label{eq::L}
L( \mathfrak X,\mu) \ = \ \mathfrak R(\mathfrak X) + \tau \cdot \mathfrak D^\circ(\mathfrak X) \; . 
\end{align}
Here, $\mathfrak R$ denotes a lower semi-continuous extrinsic risk measure and $\mathfrak D^\circ$ a lower semi-continuous intrinsic information measure. For general nonlinear systems, the parameter $\tau > 0$ can be used to tradeoff between risk and deviation. In the context of constrained linear systems, however, such a tradeoff is superfluous, as formally proven in the sections below.

\begin{remark}
The stage cost function~\eqref{eq::L} can be augmented by a control penalty. For example, one could set
\begin{align}
\label{eq::L2}
L( \mathfrak X,\mu) \ = \ \mathfrak R(\mathfrak X) + \tau \cdot \mathfrak D^\circ(\mathfrak X) + \mathfrak C(\mu) \; ,
\end{align}
where $\mathfrak C: \mathcal U \to \mathbb R$ models a lower semi-continuous control cost. This additional term does, however, not change the fact that the parameter $\tau$ does not affect the optimal solution of~\eqref{eq::dcp}. Details about how to construct $\mathfrak C$ in practice will be discussed later on in this paper, see Section~\ref{sec::Polytopes}.
\end{remark}

\subsection{Separation of Meta-Learning and Robust Control}
The goal of this section is to show that one can break the dual control problem~\eqref{eq::dcp} into an intrinsic meta learning problem and an extrinsic robust control problem. We assume that
\begin{enumerate}

\item the stage cost function $L$ has the form~\eqref{eq::L},

\item the function $\mathfrak R$ is an extrinsic risk measure,

\item the function $\mathfrak D^\circ$ is intrinsic and $\tau \geq 0$, and

\item the function $M$ is extrinsic and monotonous,
\[
\mathfrak X \subseteq \mathfrak Y \qquad \Longrightarrow \qquad M(\mathfrak X) \leq M(\mathfrak Y) \; .
\]
\end{enumerate}

\bigskip
\noindent
In this context, the meta learning problem consists of computing a constant information tube that is found by evaluating the recursion
\begin{eqnarray}
\label{eq::YRecursion}
\begin{array}{rcl}
\forall k \in \mathbb N, \qquad \mathfrak Y_{k+1} &\ \defeq \ & \overline F(\mathfrak Y_k,\nu_k) \\
\text{with} \qquad \qquad \mathfrak Y_0 & \defeq & \{ X \in \mathbb K_\mathrm{c}^n \mid X \subseteq X_0 \},
\end{array}
\end{eqnarray}
for a constant sequence $\nu_0,\nu_1, \ldots \in \mathcal U$. For simplicity of presentation, we assume $0 \in \mathbb U$ such that we can set $\nu_k(X) = 0$ without loss of generality. Due to Theorem~\ref{thm::1}, $L$  satisfies
\[
L(\mathfrak X_k) \ = \ \mathfrak R(\mathfrak X_k) + \tau \cdot \mathfrak D^\circ( \mathfrak Y_k ) 
\]
along any optimal tube of~\eqref{eq::dcp}. Consequently,~\eqref{eq::dcp} reduces to a robust control problem in the sense that all objective and constraint functions are extrinsic, while the shapes, sizes and orientations of the sets of the optimal information tube are constants, given by~\eqref{eq::YRecursion}.

\bigskip
\noindent
In summary, the contribution of intrinsic information to the objective value of~\eqref{eq::dcp}, denoted by
\[
J_\mathrm{I}(X_0) \ \defeq \ \tau \cdot \sum_{k=0}^{N-1} \mathfrak D^\circ(\mathfrak Y_k),
\]
depends on $X_0$ but it does not depend on the choice of the control law. It can be separated from the contribution of extrinsic objective terms, as elaborated below.

\subsection{Existence of Solutions}
\label{sec::existence}
In order to discuss how one can---after evaluating the meta-learning recursion~\eqref{eq::YRecursion}---rewrite~\eqref{eq::dcp} in the form of an extrinsic robust control problem, a change of variables is introduced. Let $\mathcal B_k$ denote the set of bounded functions of the form $c_k: \mathfrak Y_k \to \mathbb R^{n}$. It is a Banach space with respect to its supremum norm
\[
\Vert c_k \Vert \ \defeq \ \sup_{X \in \mathfrak Y_k} \| c_k(X) \| .
\]
Due to Theorem~\ref{thm::1}, any tight information tube $\mathfrak X_0,\mathfrak X_1, \ldots \in \mathbb I_\mathrm{c}^n$, started at $\mathfrak X_0 = \mathfrak Y_0$, is intrinsically equivalent to the precomputed tube $\mathfrak Y_0, \mathfrak Y_1, \ldots \in \mathbb I_\mathrm{c}^n$ and can be written in the form
\begin{align}
\label{eq::JC}
\mathfrak X_k = \left\{ \ Y + c_k(Y) \ \middle| \ Y \in \mathfrak Y_k \ \right\}
\end{align}
for suitable translation functions $c_k \in \mathcal B_k$. In the following, we introduce the auxiliary set
\[
\mathcal C_k \ \defeq \ \left\{ \ (c,c^+) \ \middle| \
\begin{array}{l}
\forall X \in \partial \left[ \mathfrak Y_k \cap \mathfrak V \right],\\[0.16cm] 
A c(X) - c^+( AX + \mathbb W) \in (-B \mathbb U)
\end{array}
\right\}
\]
recalling that $\partial$ denotes the boundary operator that returns the set of extreme sets of a given information ensemble. Because $\mathbb U$ is compact, $\mathcal C_k \subseteq \mathcal B_k \times \mathcal B_{k+1}$ is a closed set. Additionally, we introduce the shorthands
\begin{align}
\mathcal R_k(c_k) & \ \defeq \ \mathfrak R\left( \ \left\{ Y + c_k(Y) \ \middle| \ Y \in \mathfrak Y_k \right\} \ \right) \notag \\
\text{and} \quad \mathcal R_N(c_N) & \ \defeq \ M\left( \ \left\{ Y + c_N(Y) \ \middle| \ Y \in \mathfrak Y_N \right\} \ \right) \; . \notag
\end{align}
Since we assume that $\mathfrak R$ and $M$ are lower-semicontinuous on $\mathbb I_\mathrm{c}^n$, the functions $\mathcal R_k: \mathcal B_k \to \mathbb R$ are lower semi-continuous on the Banach spaces $\mathcal B_k$. They can be used to formulate the extrinsic robust control problem\footnote{If the sets $\mathbb U$ and $\mathbb X$ and the functions $\mathcal R_k$ are convex,~\eqref{eq::dcp2} is a convex optimization problem.}
\begin{align}
\hspace{-0.1cm} J_{\mathrm{E}}(X_0) \ = & \min_{c_0,c_1,\ldots,c_N} \; \; \sum_{k=0}^{N-1} \mathcal R_k(c_k) + \mathcal R_N(c_N)
\notag \\[0.16cm]
\label{eq::dcp2}
& \quad \ \ \mathrm{s.t.} \ \ \left\{
\begin{array}{l}
\forall k \in \{ 0, 1, \ldots, N-1 \}, \\
(c_k,c_{k+1}) \in \mathcal C_k, \ c_0 \equiv 0, \\
\forall Y \in \mathfrak Y_k, \ \ Y + c_k(Y) \subseteq \mathbb X,
\end{array}
\right.
\end{align}
which can be used to find the optimal solution of~\eqref{eq::dcp}. In detail, this result can be summarized as follows.

\bigskip
\noindent
\begin{theorem}
\label{thm::existence}
Let $\mathbb X \in \mathbb K^n$ be a closed set, let $\mathbb U \in \mathbb K_\mathrm{c}^{n_u}$, $\mathbb V \in \mathbb K_\mathrm{c}^{n_v}$, and $\mathbb W \in \mathbb K_\mathrm{c}^{n_w}$ be compact sets, let $L$ be given by~\eqref{eq::L} with $\mathfrak R$ and $M$ being set-monotonous and lower semi-continuous extrinsic risk measures, and let $\mathfrak D^\circ$ be an intrinsic lower semi-continuous information measure. Then the following statements hold.
\begin{enumerate}

\item Problem~\eqref{eq::dcp2} admits a minimizer or is infeasible.

\item Problem~\eqref{eq::dcp} is intrinsically separable; that is,
\[
J(X_0) \ = \ J_\mathrm{E}(X_0) + J_\mathrm{I}(X_0).
\]

\item If $c_0,c_1,\ldots,c_N$ is a minimizer of~\eqref{eq::dcp2}, its associated sequence of information ensembles, given by~\eqref{eq::JC}, is an optimal information tube of~\eqref{eq::dcp}.
\end{enumerate}
\end{theorem}

\bigskip
\noindent
\proof Because the objective functions $\mathcal R_k$ of~\eqref{eq::dcp2} are lower semicontinuous and since the feasible set of~\eqref{eq::dcp2} is closed under the listed assumptions, it follows directly from Weierstrass' theorem that this optimization problem admits a minimizer or is infeasible. Next, a relation between~\eqref{eq::dcp} and~\eqref{eq::dcp2} needs to be established. For this aim, we divide the proof into three parts.

\bigskip
\noindent
\textit{Part I.} Let us assume that $\mathfrak X_0,\mathfrak X_1, \ldots, \mathfrak X_N \in \mathbb I_\mathrm{c}^n$ is a tight information tube for given $\mu_0,\mu_1,\ldots,\mu_{N-1} \in \mathcal U$,
\[
\forall k \in \{ 0,1,\ldots,N-1\}, \qquad \mathfrak X_{k+1} \ = \ \overline F(\mathfrak X_k,\mu_k) \; .
\]
Due to Theorem~\ref{thm::1}, there exist functions $c_k: \mathcal B_k \to \mathbb R^n$ such that $\mathfrak X_k = \{ Y + c_k(Y) \mid Y \in \mathfrak Y_k \}$. The goal of the first part of this proof is to show that $(c_k,c_{k+1}) \in \mathcal C_k$. Because the information tube is tight, we have
\[
A X + B \mu_k(X) + \mathbb W \ \in \ \partial \mathfrak X_{k+1}
\]
for all $X \in \partial [ \mathfrak X_{k} \cap \mathfrak V ]$. Since any set $Y \in \partial [ \mathfrak Y_k \cap \mathfrak V]$ is mapped to an extreme set
\[
X = Y + c_k(Y) \in \partial [ \mathfrak X_k \cap \mathfrak V],
\]
it follows that
\begin{align}
& A (Y+c_k(Y)) + B \mu_k(X) + \mathbb W \in \partial \mathfrak X_{k+1} \notag \\[0.2cm]
\Longrightarrow \qquad & \underbrace{(AY+\mathbb W)}_{\in \ \partial \mathfrak Y_{k+1}} + \underbrace{(c_k(Y) + B \mu_k(X))}_{\in \ \mathbb R^n} \ \in \ \partial \mathfrak X_{k+1} \notag
\end{align}
for any such pair $(X,Y)$.
But this is only possible if
\begin{align}
\label{eq::cincl}
c_k(Y) + B  \mu_k(X) \ = \ c_{k+1}(AY+\mathbb W) \; .
\end{align}
Since $\mu_k(X) \in \mathbb U$ and since the choice of $Y \in \partial [ \mathfrak Y_k \cap \mathfrak V]$ is arbitrary, it follows from~\eqref{eq::cincl} that $(c_k,c_{k+1}) \in \mathcal C_k$.

\bigskip
\noindent
\textit{Part II.} The goal of the second part of this proof is to reverse the construction from the first part. For this aim, we assume that we have functions $c_k: \mathcal B_k \to \mathbb R^n$ that satisfy the recursivity condition $(c_k,c_{k+1}) \in \mathcal C_k$ for all $k \in \{ 0,1,\ldots,N-1\}$ while the sets $\mathfrak X_k$ are given by~\eqref{eq::JC}. Since every set $X \in \mathfrak X_k \cap \mathfrak V$ is contained in at least one extreme set $\overline X \in \partial \left[ \mathfrak X_k \cap \mathfrak V \right]$, there exists for every such $X$ a set $\overline Y \in \partial [ \mathfrak Y_k \cap \mathfrak V ]$ with
\[
X \ \subseteq \ \overline X \ = \ \overline Y + c_k(\overline Y).
\]
Note that this is equivalent to stating that there exists a function $\Sigma_k: \mathfrak X_k \cap \mathfrak V \to \partial [ \mathfrak Y_k \cap \mathfrak V ]$ that satisfies
\begin{align}
\label{eq::AUX1}
\forall X \in \mathfrak X_k \cap \mathfrak V, \quad X \subseteq \Sigma_k(X) + c_k(\Sigma_k(X)) \; .
\end{align}
It can be used to define the control laws
\begin{align}
\label{eq::AUX2}
\hspace{-0.1cm} \mu_k(X)  \defeq B^\dagger[ c_{k+1}(A \Sigma_k(X)+\mathbb W ) - A c_k(\Sigma_k(X)) ],
\end{align}
where $B^\dagger$ denotes the pseudo-inverse of $B$. Because we assume $(c_k,c_{k+1}) \in \mathcal C_k$, we have $\mu_k(X) \in \mathbb U$ and
\begin{align}
& A X + B \mu_k(X) + \mathbb W \notag \\[0.16cm]
& \begin{array}{rcl}
& \overset{\eqref{eq::AUX1},\eqref{eq::AUX2}}{\subseteq} & A \Sigma_k(X) + \mathbb W + c_{k+1}(A \Sigma_k(X) + \mathbb W) \in \ \mathfrak X_{k+1} \notag
\end{array}
\end{align}
for all $X \in \mathfrak X_k \cap \mathfrak V$, where the latter inclusion follows from~\eqref{eq::JC} and the fact that $A \Sigma_k(X) + \mathbb W \in \mathfrak Y_{k+1}$. Consequently, we have $\mathfrak X_{k+1} \supseteq F(\mathfrak X_k,\mu_k)$. 

\bigskip
\noindent
\textit{Part III.} The construction from Part I can be used to start with any feasible information tube of~\eqref{eq::dcp} to construct a feasible sequence $c_0,c_1,\ldots,c_N$ such that
\begin{align}
J_\mathrm{E}(X_0) & \ \leq \ \sum_{k=0}^{N-1}  \mathcal R_k(c_k) + \mathcal R_N(c_N) \notag \\[0.16cm]
& \ = \ \sum_{k=0}^{N-1} L(\mathfrak X_k,\mu_k) + M(\mathfrak X_N) - J_\mathrm{I}(X_0) \; .
\end{align}
Thus, we have $J_\mathrm{E}(X_0) + J_\mathrm{I}(X_0) \leq J(X_0)$. Similarly, the construction from Part II can be used to start with an optimal solution of~\eqref{eq::dcp2} to construct a feasible point of~\eqref{eq::dcp}, which implies $J_\mathrm{E}(X_0) + J_\mathrm{I}(X_0) \geq J(X_0)$. Thus, the second and the third statement of the theorem hold.
\qed

\subsection{Recursive Feasibility and Stability}
Feasible invariant information ensembles $\mathfrak X_\mathrm{s} \in \mathbb I_\mathrm{c}^n$ exist if and only if the optimization problem
\begin{align}
\min_{\mathfrak X_\mathrm{s},\mu_\mathrm{s}} \; \; & L(\mathfrak X_\mathrm{s},\mu_\mathrm{s}) \quad
\label{eq::Invariant}
\mathrm{s.t.} \quad \left\{
\begin{array}{l}
F(\mathfrak X_\mathrm{s},\mu_\mathrm{s}) \subseteq \mathfrak X_\mathrm{s}, \\
\mu_\mathrm{s} \in \mathcal U, \\
\forall X \in \mathfrak X_\mathrm{s},  \ \ X \subseteq \mathbb X 
\end{array}
\right.
\end{align}
is feasible. By solving this optimization problem, one can find optimal invariant information ensembles avoiding the constructions from the proof of Lemma~\ref{lem::InvariantSet}; see Remark~\ref{rem::InvariantSet}. In analogy to terminal regions in traditional MPC formulations~\cite{Rawlings2017} invariant information ensembles can be used as a terminal constraint,
\[
\mathfrak X_N \subseteq \mathfrak X_\mathrm{s}.
\]
If~\eqref{eq::dcp} is augmented by such a terminal constraint, recursive feasibility can be guaranteed. Similarly, if one chooses the terminal cost $M$ such that
\[
\min_{\mu \in \mathcal U} \ L(\mathfrak X,\mu) + M \left( \overline F(\mathfrak X, \mu) \right) \ \leq \ M(\mathfrak X)
\]
for all $\mathfrak X \in \mathbb I_\mathrm{c}^n$, the objective value of~\eqref{eq::dcp} descends along the trajectories of its associated closed-loop system. Under additional assumptions on the continuity and positive definiteness of $L$, this condition can be used as a starting point for the construction of Lyapunov functions. The details of these constructions are, however, not further elaborated at this point, as they are analogous to the construction of terminal costs for traditional Tube MPC schemes~\cite{Rawlings2017,Villanueva2020}.

\section{Polytopic Approximation Methods}
\label{sec::Polytopes}
This section discusses how to solve the dual control problem~\eqref{eq::dcp} by using a polytopic approximation method. For this aim, we assume that $\mathbb V$ and $\mathbb W$ are given convex polytopes, while $\mathbb X$ and $\mathbb U$ are convex sets.

\subsection{Configuration-Constrained Polytopes}
Polytopic computing~\cite{Fukuda2020} is the basis for many set-theoretic methods in control~\cite{Blanchini2015}. Specifically, tube model predictive control methods routinely feature parametric polytopes with frozen facet directions~\cite{Rakovic2012,Rakovic2016}. In this context, configuration-constrained polytopes are of special interest, as they admit a joint parameterization of their facets and vertices~\cite{Villanueva2022}. They are defined as follows.

Let $Y \in \mathbb R^{m \times n}$ and $G \in \mathbb R^{n_G \times m}$ be matrices that define the parametric polytope
\begin{align}
\begin{array}{rcl}
\forall y \in \mathcal G, \qquad P(y)&  \ \defeq \ & \{ x \in \mathbb R^n \mid Y x \leq y \} \notag \\[0.16cm]
\text{on} \qquad \qquad\mathcal G & \ \defeq \ & \{ y \in \mathbb R^m \mid G y \leq 0  \} ; \notag
\end{array}
\end{align}
and let $\Lambda_1,\Lambda_2,\ldots,\Lambda_{\nu} \in \mathbb R^{n \times m}$ be vertex maps, such that
\[
P(y) = \mathrm{conv}( \Lambda_1 y, \Lambda_2 y, \ldots, \Lambda_{\nu} y ) \quad \Longleftrightarrow \quad y \in \mathcal G \; ,
\]
where $\mathrm{conv}(\cdot)$ denotes the convex hull. The condition $y \in \mathcal G$ is called a configuration-constraint. It restricts the parameter domain of $P$ to a region on which both a halfspace and a vertex representation is possible. Details on how to construct the template matrix $Y$ together with the cone $\mathcal G$ and the matrices $\Lambda_i$ can be found in~\cite{Villanueva2022}.

\subsection{Polytopic Information Ensembles}
As pointed out in Section~\ref{sec::ExtemeSets}, the minimal representation of a closed information ensemble $\mathfrak X \in \mathbb I_\mathrm{c}^n$ is given by its set of extreme sets, denoted by $\partial \mathfrak X$. This motivates to discretize~\eqref{eq::dcp}, by introducing a suitable class of information ensembles, whose extreme sets are configuration-constrained polytopes. In detail,
\[
\mathfrak P(z) \ \defeq \ \left\{ \ X \in \mathbb K_\mathrm{c}^n \ \middle| \ \exists y \in \mathbb P(z): \ X \subseteq P(y) \ \right\}
\]
defines such a class of information ensembles with the polytope
\[
\mathbb P(z) \ \defeq \ \{ \ y \in \mathbb R^{m} \ \mid \ G y \leq 0, \ Z y \leq z \ \} \subseteq \mathcal G
\]
being used to parameterize convex subsets of $\mathcal G$. The choice of $Z \in \mathbb R^{l \times m}$ influences the polytopic discretization accuracy and $z \in \mathbb R^l$ denotes its associated discretization parameter. Note that $\mathfrak P(z) \in \mathbb I_\mathrm{c}^n$ is for any such $z$ a compact but potentially empty information ensemble.

\subsection{Polytopic Meta Learning}
\label{sec::intersection}
Traditional set-theoretic methods face a variety of computational difficulties upon dealing with output feedback problems, as summarized concisely in~\cite[Chapter 11]{Blanchini2015}. The goal of this and the following sections is to show that the proposed meta learning framework has the potential to overcome these difficulties. Here, the key observation is that Proposition~\ref{prop::meta} alleviates the need to intersect infinitely many information sets for the sake of predicting the evolution of a learning process. Instead, it is sufficient to compute one intersection at meta level in order pass from a prior to a posterior information ensemble.

\bigskip
\noindent
In detail, if we assume that our prior information about the system's state is represented by a polytopic ensemble, $\mathfrak P(z)$, the posterior
\[
\mathfrak P(z) \sqcap \mathfrak V \ = \ \mathfrak P(z) \cap \mathfrak V
\]
needs to be computed, where $\mathfrak V $ is given by~\eqref{eq::V}. Since $\mathbb V$ is assumed to be a polytope, $\mathfrak V$ can be written in the form
\[
\mathfrak V \ = \ \left\{
\ X \in \mathbb K_\mathrm{c}^n \ \middle| \ 
\exists y \in \mathcal G: \ X \subseteq P(y), \ Z_1 y \leq \overline v
\
\right\},
\]
as long as the template matrices $Y$, $G$, and $Z_1 \in \mathbb R^{l_1 \times m}$ as well as the vector $\overline v \in \mathbb R^{l_1}$ are appropriately chosen. Here, we construct the matrix $Z = \left(Z_1^\tr,Z_2^\tr \right)^\tr$ such that its first $l_1 \leq l$ rows coincide with $Z_1$. This particular construction of $Z$ ensures that the intersection
\[
\mathfrak P(z) \cap \mathfrak V \ = \ \mathfrak P(\zeta) \qquad \text{with} \qquad \left\{
\begin{array}{l}
\zeta_1 = \min( z_1, \overline v) \\
\zeta_2 = z_2
\end{array}
\right.
\]
can be computed explicitly, where $\min(z_1,\overline v)$ denotes the componentwise minimizer of the vectors $z_1$ and $\overline v$. The latter condition is not jointly convex in $z$ and $\zeta$. Therefore, the following constructions are based on the convex relaxation
\begin{align}
\label{eq::IntersectionRelaxation}
\left(
\begin{array}{c}
\overline v \\
z_2
\end{array}
\right) \ \leq \ \left(
\begin{array}{c}
\zeta_1 \\
\zeta_2
\end{array}
\right) \quad \Longrightarrow \quad 
\mathfrak P(z) \cap \mathfrak V \ \subseteq \ \mathfrak P(\zeta) \; .
\end{align}
Note that the conservatism that is introduced by this convex relaxation is negligible if the measurement error set $\mathbb V$ is small. In fact, for the exact output feedback case, $\mathbb V = \{ 0 \}$, we have $\min(z_1,\overline v) = \overline v$, since the measurements are exact and, as such, always informative.

\subsection{Extreme Vertex Polytopes}
In analogy to the construction of the domain $\mathcal G$, a configuration domain
\[
\mathcal H \ = \ \{ \ z \in \mathbb R^l \ \mid \ H z \leq 0 \ \}  
\]
can be chosen. In detail, by using the methods from~\cite{Villanueva2022}, a matrix $H \in \mathbb R^{l \times n_H}$ and matrices $\Omega_1,\ldots,\Omega_{\overline \nu} \in \mathbb R^{m \times l}$ can be pre-computed, such that
\[
\mathbb P(\zeta) = \mathrm{conv}( \Omega_1 \zeta, \Omega_2 \zeta, \ldots, \Omega_{\overline \nu} \zeta ) \quad \Longleftrightarrow \quad \zeta \in \mathcal H \; .
\]
This has the advantage that the vertices $\Omega_j \zeta$ of the polytope $\mathbb P(\zeta)$ are known. In order to filter the vertices that are associated to extreme polytopes, the index set
\[
\mathbb J \ \defeq \ \left\{ \ j \in \{ 1,2,\ldots,\overline \nu \} \ \middle| \ P(\Omega_j \zeta) \ \in \ \partial [\mathfrak P(\zeta)]  \ \right\}
\]
is introduced. Its definition does not depend on the choice of the parameter $\zeta \in \mathcal H$. This follows from the fact that the normal cones of the vertices of $\mathbb P(\zeta)$ do not depend on $\zeta \in \mathcal H$---recalling that the facet normals of $\mathbb P(\zeta)$ are given constants.

\bigskip
\noindent
\begin{definition}
The polytopes $P(\Omega_j \zeta)$, with $j \in \mathbb J$, are called the extreme vertex polytopes of $\mathfrak P(\zeta)$.
\end{definition}

\bigskip
\noindent
Extreme vertex polytopes play a central role in the context of designing polytopic dual control methods. This is because their shapes, sizes, and orientations can be interpreted as representatives of the intrinsic information of $\mathfrak P(\zeta)$. Moreover, the convex hull of the extreme vertex polytopes encodes the extrinsic information of $\mathfrak P(\zeta)$,
\begin{align}
\label{eq::extremeConvex}
\mathrm{conv}\left( \ \left\{ \ P(\Omega_j \zeta) \ \mid \ j \in \mathbb J \ \right\} \ \right) \ = \ \bigcup_{X \in \mathfrak P(\zeta)} X \; .
\end{align}
The latter equation follows from the fact that the vertices of the extreme polytopes of $\mathfrak P(\zeta)$ are contained in the convex hull of the vertices $\Lambda_i \Omega_j \zeta$ of the extreme vertex polytopes, with $i \in \{ 1, \ldots, \nu \}$ and $j \in \mathbb J$.

\subsection{Polytopic Information Tubes}
The goal of this section is to show that it is sufficient to assign one extreme control input $u_j \in \mathbb U$ to each extreme vertex polytope $P(\Omega_j \zeta)$ in order to discretize the control law, without introducing conservatism. This construction is similar in essence to the introduction of the vertex control inputs that are routinely used to compute robust control invariant polytopes~\cite{Gutman1986,Langson2004,Villanueva2022}. The key difference here, however, is that the ``vertices'' $P(\Omega_j \zeta)$ of the information ensemble $\mathfrak P(\zeta)$ are sets rather than vectors. They represent possible realizations of the system's information state, not a state estimate.

\bigskip
\noindent
Let us assume that $\mathbb W = P(\overline w)$ is a polytope with given parameter $\overline w \in \mathcal G$. Moreover, we assume that the vertices of $\mathbb P(\cdot)$ are enumerated in such a way that
\[
\mathbb J = \{ 1,2,\ldots, |\mathbb J|\},
\]
where $|\mathbb J| \leq \overline \nu$ denotes the cardinality of $\mathbb J$. Let us introduce the convex auxiliary set
\[
\mathcal F \ \defeq \ \left\{
(z,z^+) \middle|
\begin{array}{l}
\exists (\zeta,\xi,u) \in \mathbb R^{l} \times (\mathbb R^m)^{|\mathbb J|} \times \mathbb U^{|\mathbb J|} \\[0.16cm] \forall i \in \{ 1, \ldots, \nu \}, \forall j \in \mathbb J, \\[0.16cm]
\overline v \leq \zeta_1, \ z_2 \leq \zeta_2, \\[0.16cm]
Y A \Lambda_i \Omega_j \zeta + YB u_{j} + \overline w \leq \xi_j, \\[0.16cm]
G \xi_j \leq 0, \ H \zeta \leq 0, \ Z \xi_j \leq z^+
\end{array}
\right\}.
\]
The rationale behind the convex constraints in this definition can be summarized as follows.
\begin{itemize}
\addtolength{\itemsep}{2pt}

\item We start with the current information ensemble $\mathfrak P(z)$.

\item The constraints $\overline v \leq \zeta_1$ and $z_2 \leq \zeta_2$ subsume~\eqref{eq::IntersectionRelaxation}.

\item The constraint $H \zeta \leq 0$ ensures that the vertices of $P(\Omega_j \zeta)$ are given by $\Lambda_i \Omega_j \zeta$, with $i \in \{ 1, \ldots, \nu \}$.

\item The extreme controls $u_j$ are used to steer all vertices of $P(\Omega_j \zeta)$ into the auxiliary polytope $P(\xi_j)$.

\item And, finally, the constraints $G \xi_j \leq 0$ and $Z \xi_j \leq z^+$ ensure that $P(\xi_j)$ is contained in $\mathfrak P(z^+)$.
\end{itemize}

\bigskip
\noindent
The above points can be used as road-map for the rather technical proof of the following theorem.

\bigskip
\noindent
\begin{theorem}
\label{thm::F}
Let $\mathcal F$ and $\mathfrak V$ be defined as above, recalling that $\mathbb W = P(\overline w)$ denotes the uncertainty set and that $\mathbb U$ is assumed to be convex. Then, the implication
\[
(z,z^+) \in \mathcal F \quad \Longrightarrow \quad \exists \mu \in \mathcal U, \quad \mathfrak P(z^+) \supseteq F(\mathfrak P(z),\mu)
\]
holds for all $z,z^+ \in \mathbb R^l$.
\end{theorem}

\bigskip
\noindent
\proof
Let us assume that $(z,z^+) \in \mathcal F$. As discussed in Section~\ref{sec::intersection}, the inequalities $\overline v \leq \zeta_1$ and $z_2 \leq \zeta_2$ in the definition of $\mathcal F$ ensure that $\mathfrak P(z) \sqcap \mathfrak V \ \subseteq \ \mathfrak P(\zeta)$. Moreover, there exists for every $X \in \mathfrak P(\zeta)$ a $y \in \mathbb P(\zeta)$ with
\[
X \ \subseteq \ P(y) \ \in \ \partial [ \mathfrak P(\zeta)] \; .
\]
Next, since we enforce $H \zeta \leq 0$, $y$ is in the convex hull of the extreme vertices. That is, there exist scalar weights $\theta_1,\theta_2,\ldots, \theta_{|\mathbb J|} \in [0,1]$ with
\[
\sum_{j \in \mathbb J} \theta_j \ = \ 1 \qquad \text{and} \qquad y \ = \ \sum_{j \in \mathbb J} \theta_j \Omega_j \zeta \; ,
\]
keeping in mind that these weights depend on $X$. They can be used to define the control law
\[
\mu(X) \ \defeq \ \sum_{j \in \mathbb J} \theta_j u_{j} \ \in \ \mathbb U \quad \text{and} \quad \xi \ \defeq \ \sum_{j \in \mathbb J} \theta_j \xi_{j} \ \in \ \mathcal G
\]
where $u_1,u_2,\ldots,u_{|\mathbb J|} \in \mathbb U$ are the extreme control inputs and $\xi_1,\xi_2,\ldots,\xi_{|\mathbb J|} \in \mathcal G$ are the auxiliary variables that satisfy the constraints from the definition of $\mathcal F$.
Note that this construction is such that the vertices of the polytope $P(y)$, which are given by $\Lambda_i y$, satisfy
\begin{eqnarray}
A \Lambda_i y + B \mu(X) + w = \hspace{-0.1cm} \sum_{j \in \mathbb J} \theta_j \hspace{-0.1cm} \left[ A \Lambda_i  \Omega_j \zeta + B u_j + w \right] \in  P(\xi), \notag
\end{eqnarray}
where the latter inclusion holds for all $w \in W$. Consequently, since this holds for all vertices of $P(y)$, we have
\begin{eqnarray}
A X + B \mu(X) + \mathbb W  &\ \subseteq \ & A P(y) + B \mu(X) + \mathbb W \ \subseteq \ P(\xi)  . \notag
\end{eqnarray}
Moreover, the above definition of $\xi$ and the constraints $G \xi_j \leq 0$ and $Z \xi_j \leq z^+$ from the definition of $\mathcal F$ imply that $Z \xi \leq z^+$ and $P(\xi) \in \mathfrak P(z^+)$. But this yields
\[
F(\mathfrak P(z),\mu) \ \subseteq \ \mathfrak P(z^+) \; ,
\]
which completes the proof.
\qed

\subsection{Polytopic Dual Control}
In order to approximate the original dual control problem~\eqref{eq::dcp} with a convex optimization problem, we assume that the stage and terminal cost functions have the form
\[
L(\mathfrak P(z),\mu) = \mathfrak l(z,u) \quad \text{and} \quad M(\mathfrak P(z)) = \mathfrak m(z)
\]
for given convex functions $\mathfrak l$ and $\mathfrak m$, where the stacked vector $u = \left( u_1^\tr,u_2^\tr,\ldots,u_{|\mathbb J|} \right)^\tr$ collects the extreme control inputs. Due to Theorem~\ref{thm::F} a conservative approximation of~\eqref{eq::dcp} is given by
\begin{eqnarray}
&\min_{z,\zeta,\xi,u} \ & \; \sum^{N-1}_{k=0} \mathfrak l(z_k,u_k) + \mathfrak m(z_N) \notag \\[0.16cm]
\label{eq::qp}
&\text{s.t.} \ & \; \left\{ \; 
\begin{array}{l}
\forall k \in \{ 0, \ldots, N-1\}, \\
\forall i \in \{ 1, \ldots, \nu \}, \ \forall j \in \mathbb J, \\
\overline v \leq \zeta_{k,1}, \ z_{k,2} \leq \zeta_{k,2}, \\
Y A \Lambda_i \Omega_j \zeta_k + YB u_{k,j} + \overline w \leq \xi_{k,j}, \\
G \xi_{k,j} \leq 0, \ H \zeta_k \leq 0, \ Z \xi_{k,j} \leq z_{k+1}, \\
u_{k,j} \in \mathbb U, \ \Lambda_i \Omega_j \zeta_k \in \mathbb X, \ Z \hat y \leq z_0 \ .
\end{array}
\right.
\end{eqnarray}
Since $\mathbb U$ and $\mathbb X$ are convex sets, this is a convex optimization problem. Its optimization variables are the parameters $z_k \in \mathbb R^l$ of the polytopic information tube, the associated extreme control inputs $u_{k,j} \in \mathbb U$ and the auxiliary variables $\zeta_k \in \mathcal H$ and $\xi_k \in \mathcal G$, needed to ensure that
\[
\forall k \in \{ 0,1,\ldots, N-1 \}, \quad F(\mathfrak P(z_k),\mu_k) \ \subseteq \ \mathfrak P(z_{k+1}) \; .
\]
Here, $X_0 = P(\hat y)$ denotes the information set at the current time, modeled by the parameter $\hat y \in \mathcal G$. The constraint $Z \hat y \leq z_0$ corresponds to the initial value condition \mbox{$X_0 \in \mathfrak P(z_0)$}. Additionally, it is pointed out that the extrinsic information content of the auxiliary ensemble $\mathfrak P(\zeta) \supseteq \mathfrak P(z) \sqcap \mathfrak V$ overestimates the extrinsic information content of $\mathfrak P(z)$. Thus, the extrinsic state constraints can be enforced by using the implication chain
\begin{align}
\forall i \in \{ 1, \ldots, \nu \}, \ \forall j \in \mathbb J, \quad \Lambda_i \Omega_j \zeta_k \in \mathbb X \qquad \notag \\[0.1cm]
\Longrightarrow \quad \bigcup_{X \in \mathfrak P(\zeta)} X \ \subseteq \mathbb X \quad \Longrightarrow \quad \bigcup_{X \in \mathfrak P(z)} X \ \subseteq \mathbb X \; .
\end{align}
Finally,~\eqref{eq::qp} is solved online whenever a new information set $X_0 = P(\hat y)$ becomes available, denoting the optimal extreme controls by $u_{k,j}^\star$. A corresponding feasible control law can then be recovered by setting
\begin{align}
& \mu_0^\star(X_0) \ \defeq \ \sum_{j \in \mathbb J} \theta_j^\star(X_0) u_{0,j}, \notag 
\end{align}
where the scalar weights $\theta_j^\star(X_0)$ can, for instance, be found by solving the convex quadratic programming problem
\begin{align}
\theta^\star(X_0) \ \defeq \ \underset{\theta \geq 0}{\mathrm{argmin}} & \left( \sum_{j \in \mathbb J} \theta_j^\star u_{0,j} \right)^2 \notag \\
 \text{s.t.} \ & \left\{
\begin{array}{l}
\sum_{j \in \mathbb J} \theta_j \Omega_j \zeta_0^\star = \hat y \\
\sum_{j \in \mathbb J} \theta_j = 1,
\end{array}
\right. \notag 
\end{align}
although, clearly, other choices for the weights $\theta_j^\star$ are possible, too. Finally, the receding horizon control loop from Section~\ref{sec::RHC} can be implemented by using the above expression for $\mu_0^\star$, while the information set update and propagation step can be implemented by using standard polytopic computation routines~\cite{Blanchini2015}.

\begin{remark}
By solving the convex optimization problem
\begin{eqnarray}
&\min_{z^\mathrm{s},\zeta^\mathrm{s},\xi^\mathrm{s},u^\mathrm{s}} \ & \; \mathfrak l(z^\mathrm{s},u^\mathrm{s}) \notag \\
\label{eq::qpInvariant}
&\mathrm{s.t.} \ & \; \left\{ \; 
\begin{array}{l}
\forall i \in \{ 1, \ldots, \nu \}, \ \forall j \in \mathbb J, \\
\overline v \leq \zeta_{1}^\mathrm{s}, \ z_{2}^\mathrm{s} \leq \zeta_{2}^\mathrm{s}, \\
Y A \Lambda_i \Omega_j \zeta^\mathrm{s} + YB u_{j}^\mathrm{s} + \overline w \leq \xi_{j}^\mathrm{s}, \\
G \xi_{j}^\mathrm{s} \leq 0, \ H \zeta^\mathrm{s} \leq 0, \ Z \xi_{j}^\mathrm{s} \leq z^\mathrm{s}, \\
u_{j}^\mathrm{s} \in \mathbb U, \ \Lambda_i \Omega_j \zeta^\mathrm{s} \in \mathbb X
\end{array}
\right.
\end{eqnarray}
an optimal control invariant polytopic information ensemble can be computed.
\end{remark}

\subsection{Structure and Complexity}
Problem~\eqref{eq::qp} admits a  tradeoff between the computational complexity and the conservatism of polytopic dual control. In detail, this tradeoff can be adjusted by the choice of the following variables.
\begin{enumerate}

\item The number of facets, $m$, the number of vertices, $\nu$, and the number of configuration constraints, $n_G$, depend on our choice of $Y$ and $G$. The larger $m$ is, the more accurately we can represent the system's intrinsic information content.

\item The number of information ensemble parameters, $l$, the number of extreme vertex polytopes, $|\mathbb J|$, and the number of meta configuration constraints, $n_H$, depends on how we choose $Z$ and $H$. The larger $|\mathbb J|$ is, the more degrees of freedom we have to parameterize the optimal dual control law.
\end{enumerate}

\bigskip
\noindent
In contrast to these numbers, the number of controls, $n_u$, is given. If we assume that $\mathbb U$ and $\mathbb X$ are polyhedra with $n_{\mathbb U}$ and $n_{\mathbb X}$ facets, these number are given by the problem formulation, too. Additionally, we recall that $N$ denotes the prediction horizon of the dual controller. The number of optimization variables $n_\mathrm{opt}$ and the number of constraints, $n_\mathrm{con}$, of Problem~\eqref{eq::qp} are given by
\begin{align}
n_\mathrm{opt} &= (2N+1) l + N |\mathbb J| ( n_u + m ) \notag \\
n_\mathrm{con} &= N \left( l +  |\mathbb J| ( n_G + n_H  + l  + n_{\mathbb U} + \nu ( m + n_{\mathbb X}) \right) + l. \notag 
\end{align}
In this context, however, it should also be taken into account that the constraints of~\eqref{eq::qp} are not only sparse but also possess further structure that can be exploited via intrinsic separation. For instance, the algebraic-geometric consistency conditions
\begin{align}
\begin{array}{rcl}
G Y &=& 0, \quad \Lambda_i Y = \mathbb{1}, \\
H Z &=& 0, \quad \Omega_j Z = \mathbb{1}, \quad \text{and} \quad Z_1 Y = 0 
\end{array} \label{eq::agrelation}
\end{align}
hold for all $i \in \{1,\ldots,\nu \}$ and all $j \in \mathbb J$, which can be used to re-parameterize~\eqref{eq::qp}, if a separation of the centers and shapes of the information sets is desired.

\bigskip
\noindent
Last but not least, more conservative but computationally less demanding variants of~\eqref{eq::qp} can be derived by freezing some of its degrees of freedom. For instance, in analogy to Rigid Tube MPC~\cite{Langson2004,Rakovic2012}, one can implement a \textit{Rigid Dual MPC} controller by pre-computing a feasible point $(z^\mathrm{s},\zeta^\mathrm{s},\xi^\mathrm{s},u^\mathrm{s})$ of~\eqref{eq::qpInvariant}. Next, we set
\begin{align}
\begin{array}{rclrcl}
z_k &=& Z Y \overline x_k + z^\mathrm{s}, \quad \zeta_k &=& ZY \overline x_k + \zeta^\mathrm{s}, \\
\text{and} \quad u_{j,k} &=& \overline u_k + u_j^\mathrm{s}, \qquad
\xi_{k,j} &=& Y \overline x_{k+1} + \xi_j^\mathrm{s} \; ,
\end{array} \notag
\end{align}
where $\overline x$ and $\overline u$ denote a central state- and a central control trajectory that are optimized online, subject to $\overline x_{k+1} = A \overline x_k + B \overline u_k$. By substituting these restrictions in~\eqref{eq::qp} and by using~\eqref{eq::agrelation}, the resulting online optimization problem can be simplified and written in the form
\begin{align}
\min_{\overline x, \overline u} \ & \; \sum^{N-1}_{k=0} \overline \ell( \overline x_k, \overline u_k) + \overline m(\overline x_N) \notag \\[0.16cm]
\label{eq::RDC}
\text{s.t.} \ & \; \left\{ \; 
\begin{array}{l}
\forall k \in \{ 0, \ldots, N-1\}, \\
A \overline x_k + B \overline u_{k} = \overline x_{k+1}, \\
Z \hat y \leq ZY \overline x_0 + z^\mathrm{s}, \\
\overline u_{k} \in \overline{\mathbb U}, \ \overline x_k \in \overline{\mathbb X} .
\end{array}
\right.
\end{align}
Problem~\eqref{eq::RDC} can be regarded as a conservative approximation of~\eqref{eq::qp}. The sets
\begin{eqnarray}
\overline X & \ \defeq \ & \left\{
\ \overline x \in \mathbb R^n \ \middle| \
\begin{array}{l}
\forall i \in \{ 1, \ldots, \nu\}, \forall j \in \mathbb J, \\
\overline x + \Lambda_i \Omega_j \zeta^\mathrm{s} \in \mathbb X
\end{array} \
\right\} \notag \\
\text{and} \quad \overline {\mathbb U} & \ \defeq \ & \left\{
\ x \in \mathbb R^n \ \middle| \
\begin{array}{l}
\forall j \in \mathbb J, \quad  \ \overline u + u_j^\mathrm{s} \in \mathbb U
\end{array} \
\right\} \notag
\end{eqnarray}
take the robustness constraint margins into account, while $\overline \ell$ and $\overline m$ are found by re-parameterizing the objective function of~\eqref{eq::qp}. Problem~\eqref{eq::RDC} is a conservative dual MPC controller that has---apart from the initial value constraint---the same complexity as certainty-equivalent MPC. Its feedback law depends on the parameter $\hat y$ of the initial information set $X_0 = P(\hat y)$.

\subsection{Numerical Illustration}
We consider the constrained linear control system
\begin{align}
& A = \frac{1}{4} \left(
\begin{array}{cc}
6 & 4 \\
1 & 3
\end{array}
\right), \quad B = \left(
\begin{array}{c}
0 \\
1
\end{array}
\right), \quad C = \left(
1 \ 0
\right), \notag \\[0.16cm]
& \mathbb X = \{ x \in \mathbb R^2 \mid x_2 \geq -45 \}, \quad \mathbb U = [-55,55], \notag \\[0.16cm]
& \mathbb W = \left[ -\nicefrac{1}{2}, \nicefrac{1}{2} \right]^2 \subseteq \mathbb R^2 , \quad \mathbb V = [-1,1] \; .
\end{align}
In order to setup an associated polytopic dual controller for this system, we use the template matrices
\begin{align}
Y = \left(
\begin{array}{rr}
1 & \\[0.16cm]
1 & 1 \\[0.16cm]
  & 1 \\[0.16cm]
-1 & \\[0.16cm]
-1 & -1 \\[0.16cm]
 & -1
\end{array}
\right) \ \ \ \text{and} \ \ \ \  G = \left(
\begin{array}{rrrrrr}
-1 & 1 & -1 \\[0.16cm]
& -1 & 1 & -1 \\[0.16cm] 
& & -1 & 1 & -1 \\
& & & -1 & 1 & -1 \\[0.16cm]
-1 & & & & -1 & 1 \\[0.16cm]
1 & -1 & & & & -1
\end{array}
\right) \notag
\end{align}
setting $m = \nu = n_G = 6$. Here, $Y$ and $G$ are stored as sparse matrices: the empty spaces are filled with zeros. By using analogous notation, we set
\begin{align}
Z = \left(
\begin{array}{rrrrrr}
1 & & & 1 \\[0.16cm]
& & 1 & & & 1 \\[0.16cm]
1 & \\[0.16cm]
& 1 \\[0.16cm]
& & 1 \\[0.16cm]
& & & 1 \\[0.16cm]
& & & & 1 \\[0.16cm]
& & & & & 1
\end{array}
\right) , \ \ H = \left( \begin{array}{rrrrrrrr}
   & -1  & -1 &  1 &    &    &    &  1 \\[0.16cm]
-1 &  1  &  1 & -1 &    &    &    & -1 \\[0.16cm]
 1 &     & -1 &  1 & -1 &    &    &    \\[0.16cm]
 1 &     &  1 & -1 &    &    &    & -1 \\[0.16cm]
 1 &  1  &    & -1 &    &    & -1 &    \\[0.16cm]
   & -1  &    &    &  1 & -1 &  1 &    \\[0.16cm]
-1 &  1  &    &    & -1 &  1 & -1 &    \\[0.16cm]
 1 &     &    &    & -1 &  1 & -1 &    \\[0.16cm]
 1 &     &    &    &    & -1 &  1 & -1
\end{array} \right) , \notag
\end{align}
$l = 8$, and $n_H = 9$, which can be used to represent six dimensional meta polytopes with $6+8=14$ facets and $\overline \nu = 68$ vertices. They have up to $|\mathbb J| = 60$ extreme vertex polytopes. The first row of $Z$ corresponds to the block matrix $Z_1$. It can be used to represent the set $\mathfrak V$ by setting $\overline v = 2$, since the diameter of $\mathbb V$ is equal to $2$. Moreover, due to our choice of $Y$ and $\mathbb W$, we have
\begin{align}
\overline w \ = \ \left[ \ \nicefrac{1}{2}, \ \ 1, \ \ \nicefrac{1}{2}, \ \ \nicefrac{1}{2}, \ \ 1, \ \ \nicefrac{1}{2} \ \right]^\tr \in \mathcal G \; . \notag 
\end{align}
Next, we construct a suitable stage cost function of the form~\eqref{eq::L2}. We choose the extrinsic risk measure
\begin{align}
\mathfrak R(\mathfrak X) \defeq  \sum_{i=1}^6  \left( \max_{\{x\}\in \mathfrak X} \ Y_i x \right)^2  \hspace{-0.17cm} + 50 \cdot \sum_{i=1}^2 \max_{\{x\},\{ x' \} \in \mathfrak X} (x_i-x_i')^2 \hspace{-0.05cm} \notag
\end{align}
and the intrinsic information measure
\begin{align}
\mathfrak D^\circ(\mathfrak X) \ \defeq \  \sum_{i=1}^2  \left( \max_{X \in \mathfrak X} \ \max_{x,x' \in X} \ \left| x_i-x_i' \right| \right)^2 \; .\notag
\end{align}
This particular choice of $\mathcal R$ and $\mathcal D^\circ$ is such that
\[
\mathfrak r(z) \ \defeq \ \mathfrak R( \mathfrak P(z) ) \quad \text{and} \quad \mathfrak d^\circ(z) \ \defeq \ \mathfrak D^\circ( \mathfrak P(z) )
\]
are convex quadratic forms in $z$ that can be worked out explicitly. Namely, we find that
\begin{eqnarray}
& & \mathfrak r(z) \ = \ \left( \sum_{i=1}^6 z_{i+2}^2 \right) + 50 \cdot \left[ (z_3+z_6)^2 + (z_5+z_8)^2 \right] \notag \\[0.16cm]
& & \text{and} \quad \mathfrak d^\circ(z) = z_1^2 + z_2^2 \; . \notag
\end{eqnarray}
Last but not least a control penalty function needs to be introduced, which depends on the extreme control inputs, for instance, we can set
\[
\mathfrak c(u) = \sum_{i=1}^{|\mathbb J|} \left[ u_i^2 + 50 \cdot \left( u_i - \frac{1}{|\mathbb J|} \sum_{j=1}^{|\mathbb J|} u_j \right)^2 \, \right]
\]
in order penalize both the extreme inputs as well as the distances of these extreme inputs to their average value. The final stage cost is given by
\[
\mathfrak l(z,u) = \mathfrak r(z) + \tau \cdot \mathfrak d^\circ(z) + \mathfrak c(u),
\]
where we set the intrinsic regularization to $\tau = 0.01$.

\begin{figure*}
\begin{center}
\includegraphics[scale=0.22]{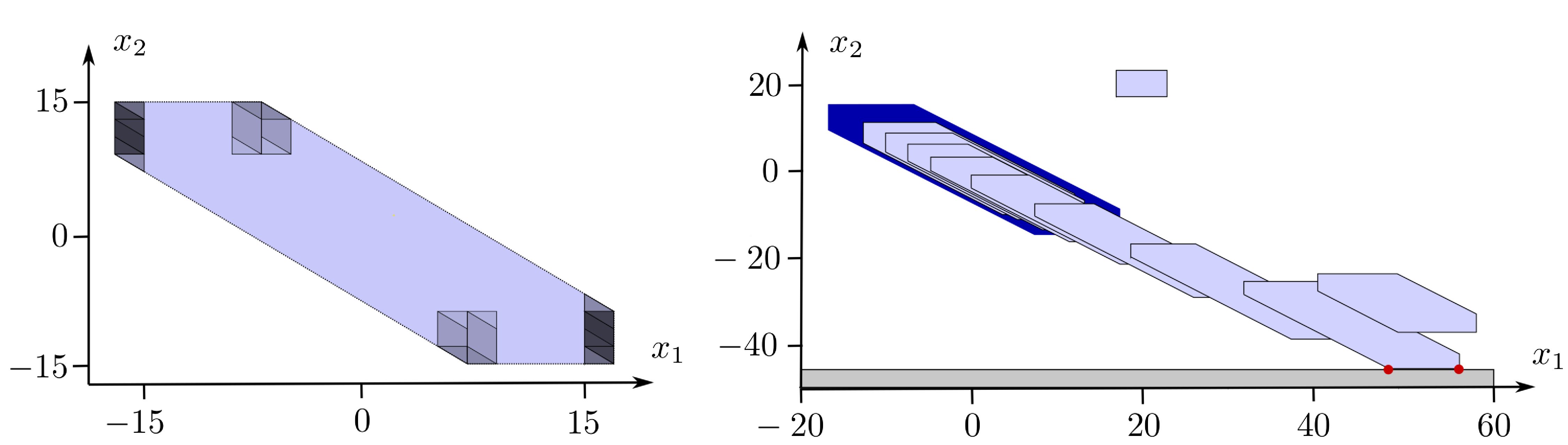}
\end{center}
\caption{\label{fig::plot} \texttt{LEFT:} Visualization of the $60$ extreme vertex polytopes of $\mathfrak P(\zeta_\mathrm{s})$, colored in gray. The convex hull of these extreme polytopes is colored in light blue. It corresponds to the union of all sets in $\mathfrak P(z_\mathrm{s})$, as predicted by~\eqref{eq::extremeConvex}.
\texttt{RIGHT:} Visualization of the extrinsic tube (light blue shaded polytopes) that corresponds to the first prediction of the dual controller~\eqref{eq::qp} with horizon $N=10$. These sets correspond to the union of the sets of the information ensembles $\mathfrak P(z_k)$. The union of the sets in the first information ensemble, $\mathfrak P(z_0)$, happens to coincide with the initial information set $X_0 = [17,23] \times [17,23]$. The extrinsic state constraints are active at the third polytope---the active vertices at which this constraint is active are colored red. The union of the sets of the optimal invariant information ensemble is shown in dark blue. It is used as terminal region.}
\end{figure*}

\bigskip
\noindent
The optimal invariant information ensemble $\mathfrak P(z_\mathrm{s})$ is found by solving~\eqref{eq::qpInvariant}. It is visualized in the left plot of Figure~\ref{fig::plot}. Note that the light blue shaded hexagon corresponds to the union of all sets in $\mathfrak P(z_\mathrm{s})$, which can be interpreted as a visualization of its extrinsic information content. The $60$ extreme vertex polytopes of $\mathfrak P(\zeta_\mathrm{s})$, given by $P(\Omega_j \zeta_\mathrm{s})$ for $j \in \{ 1,2,\ldots, 60 \}$, are difficult to plot as they are all clustered at the vertices of the extrinsic hexagon (they partly obscure each other; not all $60$ are clearly visible), but an attempt is made to visualize them in different shades of gray. As the optimal solution happens to satisfy $\mathfrak P(z_\mathrm{s})\cap \mathfrak V = \mathfrak P(\zeta_\mathrm{s})$, at least for this example, the convex relaxation~\eqref{eq::IntersectionRelaxation} does not introduce conservatism.

\bigskip
\noindent
Next, a closed-loop simulation of the polytopic dual controller~\eqref{eq::qp} is started with the initial information set
\[
X_0 = [17,23] \times [17,23]
\]
using the prediction horizon $N=10$ while the terminal cost is set to
\[
\mathfrak m(z_N) = \left\{
\begin{array}{ll}
0 & \text{if} \quad z_N \leq z_\mathrm{s} \\
\infty \ & \text{otherwise}
\end{array}
\right.
\]
in order to enforce recursive feasibility. The right plot of Figure~\ref{fig::plot} shows an extrinsic image of the first predicted tube; that is, the union of the sets $\mathfrak P(z_k)$ along the optimal solution of~\eqref{eq::qp} for the above choice of $X_0$, which are shown in light blue. The dark blue shaded polytope corresponds to the terminal region that is enforced by the above choice of $\mathfrak m$. 

\begin{remark}
The proposed polytopic dual control method optimizes feedback laws that depend on the system's information state. Note that such dual control laws are, in general, less conservative than robust output feedback laws that are based on state estimation or observer equations with affine structure, as considered in~\cite{Goulart2007} or~\cite{Mayne2009}.
\end{remark}

\section{Conclusions}
\label{sec::conclusions}

This paper has presented a set-theoretic approach to dual control. It is based on meta information spaces that enable a data-free algebraic characterization of both the present and the future information content of learning processes. In detail, an intrinsic equivalence relation has been introduced in order to separate the computation of the future information content of a constrained linear system from the computation of its robust optimal control laws. An associated intrinsic separation principle is summarized in Theorem~\ref{thm::1}. It is the basis for analyzing the existence of solutions of a large class of dual control problems under certain structural and continuity assumptions that are summarized in Theorem~\ref{thm::existence}.

\bigskip
\noindent
For the first time, this paper has presented a polytopic dual control method for constrained linear systems that is based on convex optimization. In contrast to existing robust output-feedback control schemes, this method optimizes control laws that depend on the system's information state. This alleviates the need to make control decisions based on state estimates or observer equations that induce a potentially sub-optimal feedback structure. Instead,~\eqref{eq::qp} optimizes a finite number of control inputs that are associated to the extreme vertex polytopes of the predicted information ensembles.

\bigskip
\noindent
A numerical case study for a system with two states has indicated that~\eqref{eq::qp} can be solved without numerical problems for moderately sized problems. For larger systems, however, the computational complexity of accurate dual control can become exorbitant. In anticipation of this problem, this paper has outlined strategies towards reducing the computational complexity at the cost of more conservatism. For instance, the Rigid Dual MPC problem~\eqref{eq::RDC} has essentially the same online complexity as a comparable certainty-equivalent MPC problem. The development of more systematic methods to tradeoff conservatism and computational complexity of polytopic dual control methods as well as extensions of polytopic dual control for constrained linear systems that aim at simultaneously learning their state and their system matrices $A$, $B$, and $C$ appear to be challenging and practically relevant directions for future research.

\end{document}